\documentclass[3p]{elsarticle}

\usepackage[utf8]{inputenc}
\usepackage[english]{babel}
\usepackage{amsthm}
\usepackage{amssymb}
\usepackage{mathtools}
\usepackage{amsmath}
\usepackage{diagbox}
\usepackage{caption}
\usepackage{graphicx}
\usepackage{nccmath}
\usepackage{color}
\usepackage{comment}
\usepackage{soul}
\usepackage{algcompatible}
\usepackage{algorithm}
\usepackage[noend]{algpseudocode}
\usepackage{bm}
\usepackage{optidef}
\usepackage{esvect}
\usepackage{makecell}
\usepackage{hyperref}

\newcommand{\me}{\mathbb{E}}
\newcommand{\mr}{\mathbb{R}}
\newcommand{\mq}{\mathbb{Q}}
\newcommand{\mz}{\mathbb{Z}}
\newcommand{\mn}{\mathbb{N}}

\newcommand{\ones}{\mathbf{1}}
\newcommand{\zeros}{\mathbf{0}}

\newcommand{\um}{\scalebox{0.5}[1.0]{\( - \)}}

\newcommand{\alim}{\lim_{n\uparrow \alpha }}

\newtheorem{mytheorem}{Theorem}
\newtheorem{mydefinition}{Definition}
\newtheorem{mylemma}{Lemma}
\newtheorem{myproposition}{Proposition}
\newtheorem{mycorollary}{Corollary}
\newtheorem{myaxiom}{Axiom}

\algnewcommand{\LineComment}[1]{\State \footnotesize \(/^*\) \textit{#1} \(^*/\) \normalsize}

\renewcommand{\hl}[1]{#1}
\newcommand{\mhl}[1]{#1}

\makeatletter
\def\ps@pprintTitle{%
  \let\@oddhead\@empty
  \let\@evenhead\@empty
  \def\@oddfoot{\reset@font\hfil\thepage\hfil}
  \let\@evenfoot\@oddfoot
}
\makeatother

\begin{document}

\title{A Non-Archimedean Interior Point Method for Solving Lexicographic Multi-Objective Quadratic Programming Problems}

\author[unipi]{Lorenzo Fiaschi}
\author[unipi]{Marco Cococcioni}

\address[unipi]{Department of Information Engineering, Largo Lucio Lazzarino 1 -- 56122, Pisa}

\begin{abstract}
This work presents a generalized implementation of the infeasible primal-dual Interior Point Method (IPM) achieved by the use of non-Archimedean values, i.e., infinite and infinitesimal numbers.
The extended version, called here non-Archimedean IPM (NA-IPM), is proved to converge in polynomial time to a global optimum and to be able to manage infeasibility and unboundedness transparently, i.e., without considering them as corner cases: by means of a mild embedding (addition of two variables and one constraint) NA-IPM implicitly and transparently manages their possible presence.
Moreover, the new algorithm is able to solve a wider variety of linear and quadratic optimization problems than its standard counterpart.
Among them, the lexicographic multi-objective one deserves particular attention, since NA-IPM overcomes the issues that standard techniques (such as scalarization or preemptive approach) have.
To support the theoretical properties of NA-IPM, the manuscript also shows four linear and quadratic non-Archimedean programming test cases where the effectiveness of the algorithm is verified.
This also stresses that NA-IPM is not just a mere symbolic or theoretical algorithm but actually a concrete numerical tool, paving the way for its use in real-world problems in the near future.
\end{abstract}

    



\maketitle

\section{Introduction}
\label{sec:intro}
Convex Quadratic Programming (QP) is a widely and deeply studied research topic \cite{nonlinear_opt_review}, with plenty of pervasive real world applications \cite{convex_opt_book}.
In spite of the availability of a huge amount of solving algorithms, one of them manifests a pleasant cross-model efficacy at the price of minimal changes in the implementation: the Interior Point Method (IPM) and its variants \cite{interior_review}.
Actually, it positively tackles a wide range of optimization tasks, which goes from linear programming to cone programming, passing through convex quadratic programming and semidefinite programming.

Another important and growing in interest type of optimization is lexicographic multi-objective (linear or non linear) programming, that is the task of seeking for optimal values of a set of functions (potentially constrained on the domain space) ranked lexicographically (the first objective has priority over the second, the second over the third, and so on)
\cite{CococcioniEtAlAMC2018,lexico1_on_EJOR,lexico2_on_EJOR,lexico4_on_EJOR,lexico3_on_EJOR}.
Common ways to tackle such a problem are scalarization \cite{lex_scalarization} and preemptive optimization \cite{lex_preemptive} (Chapter 17.7).
The scalarization approach transforms the multi-objective problem into a single objective one by computing a weighted sum of the functions according to their priority (the most important one takes a higher weight, the lowest important a lower weight) .
Unfortunately, the correct a priori choice of the weights is unknown, actually depending on the functions range and on the optimal solution as well.
The lack of any guarantee of correct convergence is not the only drawback: actually, the choice of too different, i.e., very high and very small weights may induce loss of precision and numerical instabilities.
On the contrary, a preemptive approach optimizes one function at time and uses its optimal value as an additional constraint for the subsequent tasks.
Even if this idea guarantees convergence to optimality, it is not drawbacks-free either.
Indeed, it requires to solve multiple optimization tasks (one for each function) with increasing complexity (addition of one constraint after each optimization is completed), which implies wasting of computational resources.
Furthermore, the nature of the problem can vary and get more complicated during time.
For instance, in a QP problem all the optimizations after the first one shall involve a convex constraint, turning them into quadratically constrained programming problems.
As a result, the optimization needs a different, more complex and less performing algorithm than the one it was supposed to be used at the beginning.

Recent literature has started to reconsider lexicographic problems, trying to solve the issues coming from both scalarization and preemptive approaches leveraging recent developments in numerical non-Archimedean computations \cite{Benci_J1,Benci_J2,Sergeyev_EMS_17}.
Some results involve game theory \cite{Fiaschi_J4,Fiaschi_J2}, evolutionary optimization \cite{Lai_et_al_PC,Lai_J2,Natural_computing}, linear programming tackled by Simplex-like methods \cite{CococcioniEtAlAMC2018,Cococcioni_20,IBigM_20}.
Actually, the key idea is to scalarize the problem adopting infinite and infinitesimal weights.
Such a choice of weights guarantees the satisfaction of lexicographic ordering during the optimization of the single-objective constrained (and non-Archimedean) function.
In a nutshell, this idea mixes the good properties of preemptive approaches (guarantee of lexicographic ordering satisfaction and convergence to the problem optimum, if any) and scalarization (transformation of the problem into a scalar one, solving it leveraging common and well established techniques). 
Thus, the idea of this work is close to the one presented in \cite{CococcioniEtAlAMC2018}, but this time using an IPM instead of a Simplex algorithm.
Once implemented, the new algorithm, called NA-IPM (Non-Archimedean Interior Point Method) will be able to solve not only lexicographic multi-objective linear programming problems but also quadratic ones and in a polynomial time.
Furthermore, by adding only two variables and one non-Archimedean constraint to the original problem, the NA-IPM is able to handle  infeasibility and unboundedness transparently, i.e., without the need to cope with them as corner cases within the implementation.
This means that ad-hoc routines (such as duality gap norm divergence) or more complex embeddings (e.g., homogeneous self-dual model) can be avoided, preventing possible slowing of computations and the need of different solving algorithms.
Lexicographic interpretation apart, NA-IPM is also able to tackle a huge variety of linear and quadratic non-Archimedean optimization problems as well, as opposed to the non-Archimedean extension of the Simplex algorithm in \cite{CococcioniEtAlAMC2018}. 
Examples of such problems are the non-Archimedean zero-sum games presented in \cite{Fiaschi_J4}.

In line with the recent literature, this work refers to those algorithms (as well as models) fully implemented (described) leveraging standard analysis with the adjective ``standard''.
This helps to better distinguish them from their non-Archimedean counterparts which, in the majority of the cases, are built over non-standard models (from here the dichotomy).
Actually, this is the case of NA-IPM which leverages the non-standard model known as Alpha Theory.
Equivalently, one could have used the word ``Archimedean'' in place of ``standard'', but it is far less common and it may have led to misunderstandings.
\hl{It is also worth to be stressed that the use of a standard or non-standard approach to model a problem does not imply that the problem itself is inherently standard or non-standard.
Rather, the theory within which a problem is modeled just defines the perimeter of tools one can use to tackle it, and set a limit to the efficacy of the chosen approach as a consequence.
This work proposes a non-standard theory to model some classes of optimization problems and proposes the NA-IPM as solving approach, believing that there are scenarios where their efficacy is higher than the one of standard techniques.}

The remaining of the work structures as follows: Section \ref{sec:ipm} and Section \ref{sec:at} provide the basic knowledge to understand this work,  the first reviews IPM and the theory behind it, while the second introduces the non-Archimedean model adopted (Benci's Alpha-Theory); Section \ref{sec:naipm} presents NA-IPM, discussing algorithm's theoretical properties (convergence, complexity, etc.), implementation issues and handling of infeasibility; Section \ref{sec:experiments} shows four numerical experiments which highlight NA-IPM effectiveness in both linear, quadratic, feasible and infeasible problems.
Finally, Section \ref{sec:conclusions} concludes the work with few summarizing lines.

\section{Reviewing the standard Interior Point Method for Quadratic Programming}
\label{sec:ipm}
Quadratic Programming (QP) \cite{numerical_opt_review} is the task of solving a problem of the form of Equation \eqref{eq:QP} or Equation \eqref{eq:QP_normal}\\
\vspace{3mm}
\begin{minipage}{.49\linewidth}
\vspace{3mm}
\begin{mini}|s|
    {}{\frac{1}{2}x^T Q x + c^Tx}{\label{eq:QP}}{}
    \addConstraint{Ax}{\le b}{}
\end{mini}
\end{minipage}
\begin{minipage}{.49\linewidth}
\vspace{3mm}
\begin{mini}|s|
    {}{\frac{1}{2}x^T Q x + c^Tx}{\label{eq:QP_normal}}{}
    \addConstraint{Ax}{= b}{}
    \addConstraint{x}{\ge0}{}
\end{mini}
\end{minipage}
where $x\in\mr^n$ is the unknown, $Q\in\mr^{n\times n}$ is symmetric, $c\in\mr^n$ together with $Q$ forms the cost function, while $A\in\mr^{m\times n}$ and $b\in\mr^m$ are the feasibility constraints, $m,\,n\in\mn$.
Whenever $Q$ is also positive semidefinite ($Q\succeq0$), the QP problem can be solved with polynomial complexity \cite{kozlov1979polynomial}; this will be the case study from now on.
A very famous algorithm to solve QP problems is the so called Interior Point method (IPM) \cite{interior_review}, in any of its many fashions: primal/primal-dual, feasible/infeasible, predictor-corrector or not. 
The success of IPM has been so wide and great that it pushed 
\cite{forsgren2002interior} to stress its spirit of unification, which has brought together areas of optimization firmly treated as disjoint.
Among the reasons of its practical success, one includes its astonishing efficiency for very large scale linear programming problems (where Simplex-like approaches generally stuck into) \cite{mehrotra_IPM}, the capability to exploit any block-matrix structure in the linear algebra operations \cite{gondzio_05}, as well as the possibility to easily and massively parallelize its computations achieving a significant speed-up (again, the same does not hold for Simplex-like approaches) \cite{gondzio_05}.
Empirical observations add one more positive feature: even if the best implementation of IPM is proved to find an $\varepsilon$-accurate solution in $\mathcal{O}(\sqrt{n}|\log\varepsilon|)$ iterations, in practice it performs much better than that, converging in almost a constant number of steps regardless the problem dimension $n$ \cite{colombo_08} ($\varepsilon$ is the optimality tolerance).

This work focuses on the predictor-corrector infeasible primal-dual IPM \cite{wright_pd_IPM}, which is broadly accepted today to be the most efficient one \cite{interior_review}.
The next subsections is devoted to resume the algorithm in question.
The reader already familiar with IPM might find it overly lengthy, but it is better to recall here the details of IPM standard implementation in order to more easily introduce its non-Archimedean extension in Section \ref{sec:naipm}.

\subsection{Description of the standard IPM algorithm}
First of all, the assumption that the problem is formulated according to \eqref{eq:QP_normal} holds.
If not, one needs to rewrite the problem by adding some slack variables, one for each constraint, without penalizing them in the cost function.
Presolving techniques \cite{presolving} for guaranteeing consistent problem formulations must be applied too.
Second, duality theory \cite{NLP_review} states that first order optimality conditions for QP problems (also known as KKT conditions) are the ones reported in Equation \eqref{eq:KKT}, where $\lambda$ are the dual variables and $s$ the dual slack variables (a slight modification which includes an additional parameter $\mu$ is provided in Equation \eqref{eq:KKT_IPM} and will be discussed later).\\
\vspace{3mm}
\begin{minipage}{.5\linewidth}
\vspace{3mm}
\begin{equation}
\begin{dcases*}
    Ax = b\\
    A^T\lambda+s-Qx=c\\
    XS\mathbf{1}=0\\
    (x,s) \ge 0
\end{dcases*}
\label{eq:KKT}
\end{equation}
\end{minipage}
\begin{minipage}{.5\linewidth}
\vspace{3mm}
\begin{equation}
\begin{dcases*}
    Ax = b\\
    A^T\lambda+s-Qx=c\\
    XS\mathbf{1}=\mu\mathbf{1}\\
    (x,s) \ge 0
\end{dcases*}
\label{eq:KKT_IPM}
\end{equation}
\end{minipage}
Third, iterative algorithms based on the Newton–Raphson method \cite{stoer2013introduction} find the solution of systems as the one in \eqref{eq:KKT} (positivity of $x$ and $s$ excluded) starting from a solution $(x^0,\,\lambda^0,\,s^0)$ and repeatedly solving the linear system in \eqref{eq:Newton} (Equation \eqref{eq:Newton_IPM} presents a slight modification with two additional parameters $\mu$ and $\sigma$ which will be introduced in the next lines), where $r_b = Ax-b$ and $r_c = A^T\lambda+s-Qx-c$, \hl{while $X$ and $S$ are the diagonal matrices obtained by $x$ and $s$ (the iteration apices are omitted for readability reasons).}\\
\vspace{3mm}
\begin{minipage}{.5\linewidth}
\vspace{3mm}
\begin{equation}
\begin{bmatrix}
    -Q & A^T & I\\
    A & 0   & 0\\
    S & 0   & X
\end{bmatrix}
\begin{bmatrix}
    \Delta x\\ \Delta \lambda \\ \Delta s
\end{bmatrix} =
\begin{bmatrix}
    -r_c \\ -r_b \\ -XS\mathbf{1}
\end{bmatrix}
\label{eq:Newton}
\end{equation}
\end{minipage}
\begin{minipage}{.5\linewidth}
\vspace{3mm}
\begin{equation}
\begin{bmatrix}
   -Q & A^T & I\\
    A & 0   & 0\\
    S & 0   & X
\end{bmatrix}
\begin{bmatrix}
    \Delta x\\ \Delta \lambda \\ \Delta s
\end{bmatrix} =
\begin{bmatrix}
    -r_c \\ -r_b \\ \sigma\mu\mathbf{1}-XS\mathbf{1}
\end{bmatrix}
\label{eq:Newton_IPM}
\end{equation}
\end{minipage}
This system is known as Newton's step equation, and once solved the approximated solutions of \eqref{eq:KKT} are updated accordingly: $(x^{k+1},\,\lambda^{k+1},\,s^{k+1})\leftarrow (x^k+\nu\Delta x,\, \lambda^k+\nu\Delta \lambda,\, s^k+\nu\Delta s)$ with $\nu\in(0,\,1]$ opportunely chosen.

Infeasible primal-dual IPM aims to find optimal primal-dual solutions $(\overline{x},\,\overline{\lambda},\,\overline{s})$ of \eqref{eq:KKT} by applying a variant of \eqref{eq:Newton} and modifying the search direction as well as the step lengths so that the non-negativity of $x$ and $s$ is strictly satisfied.
To do this, IPM is designed as a path-following algorithm, i.e., an algorithm which tries to stay as close as possible to a specific curve during the optimization \cite{gonzaga_92}.
Such a curve, also known as \emph{central path}, is uniquely identified by all the duality-gap values $\mu\in(0,\tfrac{x^{0T}s^0}{n}]$, where $(x^0,\,\lambda^0,\,s^0)$ is the primal-dual starting point.
Its role is to keep the intermediate solutions away from the feasible region's boundaries as much as possible (from here the word ``central").
The points belonging to the central path are those satisfying conditions in Equation \eqref{eq:KKT_IPM}, which is a slight modification of \eqref{eq:KKT}, for a given value of $\mu$.
The uniqueness of the curve is guaranteed whenever the feasible region admits at least one strictly feasible point, that is a primal-dual solution such that $(x,s)>0$ (the latter is a very weak assumption).
The Newton's step is now computed accordingly with Equation \eqref{eq:Newton_IPM}, where $\sigma\mu$ is the next duality gap to achieve, $\sigma\in(0,1)$; furthermore, $\sigma\mu$ uniquely identifies the primal-dual point on the central path to target.
Since Equation \eqref{eq:Newton_IPM} approximates \eqref{eq:Newton} more and more closely as $\mu$ goes to zero, if the latter converges to anything as $\mu\downarrow0$ (and it does if the feasible region is bounded), then IPM eventually converges to an optimal solution of the QP problem.

Predictor-corrector version of IPM splits the Newton's direction computation in two stages: predictor and corrector.
The predictor one solves the system in \eqref{eq:Newton}, i.e., it looks for a solution as optimal as possible, identifying it by the three predictive directions solving the Newton's step equation, namely $\Delta x_p$, $\Delta \lambda_p$, $\Delta s_p$ .
Then, it is the turn of the corrector stage which operates a restoration of some of the centrality lost during the predictor stage.
Indeed, pushing all the way towards the optimum may take the temporary solutions too close to the boundaries, affecting the quality of the convergence due to possible numerical issues.
To cope with that, the corrector stage solves the linear system in \eqref{eq:Newton_IPM} with the constant term $[0,\,\, 0,\,\, \sigma\mu\mathbf{1}-\Delta x_p \odot \Delta s_p]$ ($\odot$ indicates the Hadamard product).
In most of the cases, the corrector direction badly (but mildy) affects the optimality of the predictor solution, still guaranteeing an improvement with respect to the previous iterate.
Indeed, a good centrality of the temporary solutions is too crucial for a good and fast convergence to be ignored.
The final optimizing direction is just the sum of the ones found at the two stages.
A final remark: since the arbitrary a-priori choice of $\sigma$ may bias too much the performance, Mehrotra \cite{mehrotra_IPM} fruitfully proposed an adaptive setting of it to $\left(\tfrac{\mu^{\text{new}}}{\mu}\right)^3$, where $\mu^{\text{new}} = \tfrac{(x+\Delta x)^T(s+\Delta s)}{n}$.

There are three aspects left to be uncovered to exhaustively present the IPM algorithm used in this work.
They consists in three crucial issues of any iterative algorithm: problem infeasibility, starting point choice and algorithm termination.
Multiple ways to address infeasibility exist, two of the most common are iterates divergence \cite{kojima_93} and problem embedding \cite{convexP_embedding}.
The first one leverages some results on convergence \cite{wright_pd_IPM}, which state that in case of infeasibility the residuals cannot decrease under a positive constant, and therefore the iterates must diverge.
To identify this phenomenon, it is enough to check at each iteration whether the norm of $x$ and $s$ exceeds a certain threshold $\omega$, i.e., $\|(x,s)\|>\omega$.
On the contrary, the second one aims to embed the given problem in a larger one for which a feasible solution always exists (called augmented solution).
Once found, it shall tell the user whether the original problem admits a feasible solution (in which case it shall allow one to retrieve it from the augmented one) or the opposite (specifying also which of the polyhedra is infeasible: the primal, the dual or both).

The starting point choice is vital for the quality of the solution, or even for the convergence itself.
In literature, there exist plenty of papers proposing novel ways to initialize an IPM.
In this work, the approach used is the one by Mehrotra \cite{mehrotra_IPM}, which first seeks for the solution of the two least squares problems in Equation \eqref{eq:x0}  and \eqref{eq:s0} in order to satisfy the equality constraints, then executes some additional operations to guarantee strict positivity and well centering (even at the expenses of feasibility).
The complete procedure is reported in Algorithm \ref{alg:start} .\\
\vspace{2mm}
\begin{minipage}{.5\linewidth}
\vspace{3mm}
\begin{mini}|s|
{x}{x^Tx}{\label{eq:x0}}{}
\addConstraint{Ax}{=b}{}
\end{mini}
\end{minipage}
\begin{minipage}{.5\linewidth}
\vspace{3mm}
\begin{mini}|s|
{\lambda,s}{s^Ts}{\label{eq:s0}}{}
\addConstraint{A^T\lambda+s-Qx}{=c}{}
\end{mini}
\end{minipage}
\begin{algorithm}[ht]
\caption{Starting point computation for IPM}
\label{alg:start}
\begin{algorithmic}[1]
\Procedure{starting\_point}{$A$, $b$, $c$, $Q$}
\LineComment{solution of \eqref{eq:x0}}
\State $x = A^T(AA^T)^{-1}b$
\LineComment{solution of \eqref{eq:s0}}
\State $\lambda = (AA^T)^{-1}A(c+Qx)$
\State $s = c+Qx-A^T\lambda$
\LineComment{guarantee of positivity}
\State $\delta_x = \max(-\tfrac{3}{2}\min x_i,\,0)$
\State $x = x+\delta_x\mathbf{1}$
\State $\delta_s = \max(-\tfrac{3}{2}\min s_i,\,0)$
\State $s = s+\delta_s\mathbf{1}$
\LineComment{increment in centrality and well-balancing}
\State $\delta_x = \tfrac{1}{2}\tfrac{x^Ts}{s^T\mathbf{1}}$
\State $\delta_s = \tfrac{1}{2}\tfrac{x^Ts}{x^T\mathbf{1}}$
\State $x = x+\delta_x \mathbf{1}$
\State $s = s+\delta_s\mathbf{1}$
\State \Return $x,\, \lambda,\, s$
\EndProcedure
\end{algorithmic}
\end{algorithm}

Finally, the algorithm termination is demanded to three aggregated values, $\rho_1$, $\rho_2$ and $\rho_3$.
They describe the degree of KKT conditions \eqref{eq:KKT} satisfaction achieved by the current iteration.
As soon as all the three values are under a user-defined tolerance threshold $\varepsilon\in\mr^+$, the algorithm is considered converged.
Actually, $\rho_1$ and $\rho_2$ refer to primal and dual feasibility, respectively, while $\rho_3$ is in charge of the complementary degree of $x$ and $s$.
Their definition is
\begin{equation}
    \rho_1 = \frac{\|Ax-b\|}{1+\|b\|}, \quad\quad
    \rho_2 = \frac{\|A^T\lambda+s-Qx-c\|}{1+\|c\|}, \quad\quad
    \rho_3 = \frac{\mu}{1+|\frac{1}{2}x^TQx+c^Tx|}.
\label{eq:rho_def}
\end{equation}

A resume of the whole procedure just described can be found in Algorithm \ref{alg:std_IPM}.
Next section introduces the non-Archimedean tools exploited to implement the more general interior point algorithm at the core of this manuscript.

\begin{algorithm}[!htpb]
\caption{Predictor-corrector infeasible primal-dual IPM}
\label{alg:std_IPM}
\begin{algorithmic}[1]
\Procedure{standard\_IPM}{$A$, $b$, $c$, $Q$, $\varepsilon$, $\omega$, $\texttt{max\_it}$}
\LineComment{flag of correct termination}
\State $\texttt{flag}$ = False
\State $x$, $\lambda$, $s$ = starting\_point($A$, $b$, $c$, $Q$)
\For{$i=1,\,\ldots,\,\texttt{max\_it}$}
    \LineComment{compute residuals}
    \State $r_b = Ax-b \quad r_c = A^T\lambda+s-Qx-c \quad r_{\mu} = x \odot s$
    \LineComment{compute centrality, $n$ is the length of $x$}
    \State $\mu = \tfrac{r_{\mu}^T\mathbf{1}}{n}$
    \LineComment{compute KKT conditions satisfaction parameters}
    \State $\rho_1 = \tfrac{\|r_b\|}{1+\|b\|} \quad \rho_2 = \tfrac{\|r_c\|}{1+\|c\|} \quad \rho_3 = \tfrac{\mu}{1+|\tfrac{1}{2}x^TQx+c^Tx|}$
    \If{$\rho_1\le\varepsilon$ \text{and} $\rho_2\le\varepsilon$ \text{and} $\rho_3\le\varepsilon$}
        \LineComment{primal-dual feasible optimal solution found}
        \State \texttt{flag} = True
        \State \Return $x$, $\lambda$, $s$, \texttt{flag}
    \EndIf
    \LineComment{compute predictor directions solving \eqref{eq:KKT}}
    \State $\Delta x_p$, $\Delta \lambda_p$, $\Delta s_p = $ predict($A$, $b$, $c$, $Q$, $r_b$, $r_c$, $r_{\mu}$)
    \LineComment{compute predictor step size}
    \State $\mhl{\nu_{pp} = 0.99\min(\max\nolimits_{\overline{\nu}}\{\overline{\nu}\,|\, x+\overline{\nu}\Delta x_p \ge 0\},\,1)}$
    \State $\mhl{\nu_{pd} = 0.99\min(\max\nolimits_{\overline{\nu}}\{\overline{\nu}\,|\, s+\overline{\nu}\Delta s_p \ge 0\},\,1)}$
    \State $\nu=\min(\nu_{pp},\,\nu_{pd})$
    \LineComment{estimate $\sigma$}
    \State $\Tilde{x} = x+\nu\Delta x_p \quad \Tilde{s} = s+\nu\Delta s_p$
    \State $\mu^{\text{new}} = \tfrac{\Tilde{x}^T\Tilde{s}}{n}$
    \State $\sigma = (\tfrac{\mu^{\text{new}}}{\mu})^3$
    \LineComment{compute corrector directions solving the corresponding Newton's system}
    \State $\Delta x_c$, $\Delta \lambda_c$, $\Delta s_c = $ corrector($A$, $b$, $c$, $Q$, $\sigma\mu\mathbf{1}-\Delta x_p \odot \Delta s_p$)
    \LineComment{compute new direction}
    \State $\Delta x = \Delta x_p + \Delta x_c \quad \Delta \lambda = \Delta \lambda_p + \Delta \lambda_c \quad \Delta s = \Delta s_p + \Delta s_c$
    \LineComment{compute step size}
    \State $\mhl{\nu_{p} = 0.99\min(\max\nolimits_{\overline{\nu}}\{\overline{\nu}\,|\, x+\overline{\nu}\Delta x \ge 0\},\,1)}$
    \State $\mhl{\nu_{d} = 0.99\min(\max\nolimits_{\overline{\nu}}\{\overline{\nu}\,|\, s+\overline{\nu}\Delta s \ge 0\},\,1)}$
    \State $\nu=\min(\nu_p,\,\nu_d)$
    \LineComment{compute target primal-dual solution}
    \State $x = x+\nu\Delta x \quad \lambda = \lambda+\nu\Delta \lambda \quad s = s+\nu\Delta s$
    \If{$\|(x,s)\|>\omega$}
        \LineComment{Divergence detected}
        \State \Return $x$, $\lambda$, $s$, \texttt{flag}
    \EndIf
\EndFor
\State \Return $x$, $\lambda$, $s$, \texttt{flag}
\EndProcedure
\end{algorithmic}
\end{algorithm}

\section{Alpha Theory}
\label{sec:at}
This section aims to introduce the non-Archimedean model adopted to generalize IPM to NA-IPM: Benci's Alpha Theory \cite{BDNbook}.
Its key ingredient is the infinite number $\alpha$, whose definition can be given in multiple ways, as usual in Mathematics.
In this work, the axiomatization of Alpha Theory given in \cite{Benci_J2} is used.
The reason of this choice is that it provides a very plain and application-oriented introduction to the topic, which has also the property to stress the numerical nature and immediate applicability of such a theory to Data Science problems.
Then, the last part of the section is devoted to show the strict relation between standard lexicographic QP problems and a particular category of non-Archimedean QP ones.
This fact suggests a straightforward application of NA-IPM to lexicographic problems, as testified by the experimental section.

\subsection{Background in brief}
The minimal ground for Alpha Theory consists in three axioms, the first of which introduces a wider field than $\mr$, namely $\me$, which actually contains it.
\begin{myaxiom}
Exists an ordered field $\mathbb{E} \supset \mathbb{R}$ whose numbers are called Euclidean numbers.
\label{ax:existence}
\end{myaxiom}
The second one aims at describing better some of the elements in $\me$.
Actually, it states that $\me$ contains numbers which are finite (since it contains $\mr$), infinite and infinitesimal according to the following definition:
\begin{mydefinition}{} Given $\xi\in\me$, then
\begin{itemize}
    \item $\xi$ is infinite $\Longleftrightarrow$ $\forall\,n\in\mn, |\xi|>n$ 
    \item $\xi$ is finite $\Longleftrightarrow$ $\exists\, n\in\mn,$ $\tfrac{1}{n}<|\xi|<n$ 
    \item $\xi$ is infinitesimal  $\Longleftrightarrow \forall\,n\in\mn,$ $|\xi| < \tfrac{1}{n}$
\end{itemize}
\label{def:inf_fin_inf}
\end{mydefinition}
\noindent To do this, it introduces the previously mentioned value $\alpha$ by means of the numerosity function $\mathfrak{num}$, which in some sense can be seen as a particular counting function.
For the sake of rigorousness, let $ V(\mathbb{N})$ denote the superstructure on $\mathbb{N}$, namely
\[
V(\mathbb{N})=\bigcup_{i=0}^\infty V_i(\mathbb{N}),
\]
where ($\mathcal{P}$ is the power set)
\[
V_0(\mathbb{N})=\mathbb{N}, \qquad V_{i+1}(\mathbb{N})=V_i(\mathbb{N})\cup\mathcal{P}(V_i(\mathbb{N}));
\]
then set
\[
\mathcal{U}:= \{ X\in V(\mathbb{N}) \ | X \text{ is countable} \}.
\]
\begin{myaxiom}
\label{ax:numerosity}
Exists a function $\mathfrak{num}$, $\mathfrak{num}: \mathcal{U} \rightarrow \mathbb{E}$ which satisfies
     \begin{itemize}
        \item if $A$ is finite $\mathfrak{num}(A)=|A| \; $ \textnormal{(}$|\cdot|$ \textnormal{denotes the cardinality of a set)}
        \item $\mathfrak{num}(A)<\mathfrak{num}(B)$ if $A\subset B$
        \item $\mathfrak{num}(A\cup B) = \mathfrak{num}(A) + \mathfrak{num}(B) - \mathfrak{num}(A \cap B)$
        \item $\mathfrak{num}(A \times B) = \mathfrak{num}(A) \cdot \mathfrak{num}(B)$
        \item $\alpha = \mathfrak{num}(\mathbb{N})$
    \end{itemize}
\end{myaxiom}
\noindent More precisely, Axiom \ref{ax:numerosity} not only introduces the infinite number $\alpha$, but also states that its manipulation by means of algebraic functions generates other Euclidean numbers.
For instance, the following equivalences hold true in $\me$ and each one individuates a specific number:
\begin{equation*}
     \alpha \cdot (\alpha + 2) = \alpha^{2} + 2\alpha \quad \quad 0 < \frac{1}{\alpha} = \alpha^{\um 1} < \alpha^0 = 1 < \alpha^1 = \alpha < (\alpha + 1)
\end{equation*}
\begin{equation*}
   \frac{-10.0\alpha^2 +16.0 +42.0\eta^2} {5.0\alpha^2+7.0} =  -2.0 +6.0\eta^2 
\end{equation*}
where there is the convention to indicate $\alpha^{\um1}$ with $\eta$, i.e., $\eta \coloneqq \alpha^{\um1}$.

The last axiom introduces the notion of $\alpha$-limit by means of which any function $f\colon\mr^m\rightarrow\mr$ can be extended to a function $f^*\colon\me^m\rightarrow\me$ satisfying all the first order properties of the former.
This is a crucial tool of any non-standard model since it guarantees the transfer of such properties from $f$ to $f^*$.
Because of this, it is commonly referred to as Transfer Principle \cite{keisler_foundations}.
\begin{myaxiom}
\label{ax:transfer}
Every sequence $\varphi :\mn\rightarrow \mr^m$ has a unique $\alpha$-limit denoted by $\alim\varphi (n)$ which satisfies the following properties:

\begin{enumerate}
\item if $\xi \in \me^m$, then \textit{there exists a sequence }$\varphi :\mn\rightarrow \mr^m$ such that
\[
\xi =\lim_{n\uparrow \alpha }\varphi (n)
\]
\item if $\varphi (n)=n,$ then
\[
\lim_{n\uparrow \alpha }\varphi (n)=\alpha
\]
\item \textit{if } $\exists n_{0}\in \mathbb{N}$ such that $\forall n\geq n_{0},\ \varphi (n)\geq \psi (n)$, then 
\[
 \lim_{n\uparrow \alpha } \varphi (n)\geq \lim_{n\uparrow \alpha }\psi (n)
\]

\item any sequence $\varphi,\psi$
\begin{equation*}
\begin{split}
\lim_{n\uparrow \alpha }\varphi (n)+\lim_{n\uparrow \alpha }\psi (n)
&=\lim_{n\uparrow \alpha }\left( \varphi (n)+\psi (n)\right) , \\
\lim_{n\uparrow \alpha }\varphi (n)\cdot \lim_{n\uparrow \alpha }\psi (n)
&=\lim_{n\uparrow \alpha }\left( \varphi (n)\cdot \psi (n)\right).
\end{split}
\end{equation*}
\end{enumerate}
\end{myaxiom}
\noindent Leveraging Axiom \ref{ax:transfer}, any real-valued function $f:\mr^m\rightarrow \mr$ can be extended to a function $f^*:\me^m\rightarrow \me$ by setting \cite{Benci_J2}
\begin{equation}
f^*(\xi) = f^{\ast }\left( \lim_{n\uparrow \alpha }\varphi (n)\right) =\lim_{n\uparrow
\alpha }f\left( \varphi (n)\right).
\label{eq:function_extension}
\end{equation}
A similar extension is possible for multivalued functions and sets too.
In the first case, for any $f\colon\mr^m\rightarrow\mr^q$ with $f(x)=\{f_1(x),\,\ldots,\,f_q(x)\}$ and $f_i\colon\mr^m\rightarrow\mr$ $\forall i=1,\ldots,q$, $f^*\colon\me^m\rightarrow\me^q$ is defined as $f^*(\xi)=\{f_1^*(\xi),\,\ldots,\,f_q^*(\xi)\}$ where $f_i^* = \alim f_i $ $\forall i=1,\ldots,q$.
In the second case, given a set $F$, $F^*$ is defined as $F^*=\{\alim f\,|\,f\in F\}$.

Notice that the use of the word ``extension'' is fully justified by the fact that $f^*(x) = f(x)$ $\forall x\in\mr$.
In the remaining of the work, when no ambiguity is possible, the ``$\ast$" will be omitted and therefore $f$ and $f^{\ast }$ will be denoted by the same symbol.

Another important implication of Axioms 1-3 is that any $\xi\in\me\setminus\{0\}$ admits the following representation \cite{BDNbook}:
\begin{equation*} \xi=\sum_{i=1}^{\infty} \psi_i \alpha^{g_{\xi}(i)},\,\psi_1\neq0,\,i\in\mn
\label{eq:poor_e}
\end{equation*}
where $\psi_i\in\mr$, $g_{\xi}\colon\mn\rightarrow\mr$ is a monotonic decreasing function.
Finally, some definitions which are useful in the remaining of the work follow.
\begin{mydefinition}[Monosemium]
$\xi\in\me$ is called monosemium if and only if $\xi=\psi\alpha^i$, $\psi\in\mr$, $i\in\mr$.
\label{def:monosemium}
\end{mydefinition}

\begin{mydefinition}[Leading monosemium]
Given $\xi\in\me$, $\psi_1\alpha^{g_{\xi}(1)}$ is called leading monosemium.
\end{mydefinition}

\begin{mydefinition}[Leading monosemium function]
The leading monosemium function $\mathtt{lead\_mon}\colon\me\rightarrow\me$ is a non-Archimedean function which maps each Euclidean number in its leading monosemium.
\end{mydefinition}

\begin{mydefinition}[Order of magnitude]
The order of magnitude is a function $\mathcal{O}\colon\me\setminus\{0\}\rightarrow\me\setminus\{0\}$ such that $\forall\xi\in\me\setminus\{0\}$ $\mathcal{O}(\xi) = \alpha^{g_{\xi}(1)}$.
\end{mydefinition}

\begin{mydefinition}[Smallest order of magnitude]
Let one indicate with $\me_{lim}\subset\me$ the set of Euclidean numbers represented by a limited number of monosemia, i.e., $\me_{lim}\coloneqq\{\xi\in\me\setminus\{0\}\,|\,\xi=\sum_{i=1}^{l} \psi_i \alpha^{g_{\xi}(i)},\,\psi_1\neq0,\,l\in\mn\}$.
The smallest order of magnitude is a function $\mathit{o}\colon\me_{lim}\rightarrow\me_{lim}$ such that $\forall\xi\in\me_{lim}$ $\mathit{o}(\xi)=\alpha^{g_{\xi}(l)}$.
\end{mydefinition}

\noindent A common generalization of the last two functions to vectorial spaces is the following:
\begin{equation*}
    \forall \xi \in\me^n\setminus\{\zeros\}\;\; \mathcal{O}(\xi)=\max_{i=1,\,\ldots,\,n}\mathcal{O}(\xi_i)
    \quad \quad
    \forall \xi \in\me^n_{lim}\;\; \mathit{o}(\xi)=\min_{i=1,\,\ldots,\,n}\mathit{o}(\xi_i).
\end{equation*}

\begin{mydefinition}[Infinitesimal number]
A way to indicate that $\xi\in\me$ is infinitesimal (equivalent to the one in Definition \ref{def:inf_fin_inf}) is by the inequality $\mathcal{O}(\xi)<\alpha^0$, which is usually indicated as $\xi\approx0$.
\end{mydefinition}

\subsection{How to handle lexicographic optimization}
This subsection aims to introduce Proposition \ref{prop:convertion}, which builds a bridge between standard lexicographic QP problems and non-Archimedean quadratic optimization.
First, let one formally introduce standard and non-Archimedean optimization models.
Then, the result comes straightforwardly.
\begin{mydefinition}[Standard optimization problem]
An optimization problem of the form
\begin{mini}|s|
{}{f(x)}{\label{eq:opt}}{}
\addConstraint{g_i(x)}{\le0}{\quad i=1,\,\ldots,\,n}
\addConstraint{h_j(x)}{=0}{\quad j=1,\,\ldots,\,m}
\end{mini}
\end{mydefinition}
is said to be standard (or Archimedean) if and only if $f$, and all $g_i$ and $h_j$ are real-valued multivariate functions.

\begin{mydefinition}[Non-Archimedean optimization problem]
An optimization problem of the form \eqref{eq:opt}
is said to be non-Archimedean if and only if at least one of the functions $f$, $g_i$ or $h_j$ is non-Archimedean, i.e., at least one of them has values in $\me^k$ (with $k\in\mn$) or maps the input into values in $\me$.

In the case of QP problems as \eqref{eq:QP} and \eqref{eq:QP_normal}, this is equivalent to say that either $x\in\me^k$ or at least one entry of $A$, $b$, $c$ or $Q$ is non-Archimedean.
\end{mydefinition}

\begin{myproposition}
\label{prop:convertion}
Consider a lexicographic optimization problem whose objective functions $f_1,\,\ldots,\,f_n$ are real functions and the priority is induced by the natural order.
Then, there exists an equivalent scalar program over the same domain, whose objective function is non-standard and has the following form: 
\[
F(x) = \beta_1 f_1(x) + \ldots + \beta_n f_n(x),
\]
where $\frac{\beta_{i+1}}{\beta_i}\approx 0$, $i=1,\,\ldots,\,n$-1.
\begin{proof}
The proof will show that the two programs share the same global maxima. Very similar considerations hold for local maxima, global and local minima, but are omitted for brevity.

Let $\Omega$ be the domain of the two programs, and $\omega\in\Omega$ be a global maximum of the lexicographic optimization problem, i.e., $\nexists \omega'\in\Omega$ such that $f_1(\omega')>f_1(\omega)$, or $\nexists i=2,\,\ldots,\,n$ such that $f_j(\omega')=f_j(\omega)$, $j=1,\,\ldots,\,i$-1, and $f_i(\omega')>f_i(\omega)$.
However, this is true if and only if $\nexists \omega'\in\Omega$ such that $F(\omega')>F(\omega)$, by the definition of $F$ itself.
The latter fact exactly means that $\omega$ is a global maximum for the scalar non-standard program.
\end{proof}
\end{myproposition}

The choice of weights $\{\beta_i\}_{i=1}^n$ adopted in this work is $\beta_i = \alpha^{1-i}$, $i=1,\,\ldots,\,n$.
As an example, consider the two-objective lexicographic optimization problem in \eqref{eq:lex_example}.
It seeks for the point in the unitary cube which minimizes the first objective $f_1(x)=x_1^2+x_2^2$.
In case more than one such point exist, as in this example where the whole segment $\mathcal{S}=\{x\in\mr^3\,|\,x_1=x_2=0,\,x_3\in[0,1]\}$ optimizes $f_1$, then such set of candidate optimal solutions is refined considering the second objective $f_2(x)=x_3+x_2$.
It is easy to verify that $f_2$ selects as optimal only one point in $\mathcal{S}$, namely the origin.

According to Proposition \ref{prop:convertion}, a possible non-Archimedean reformulation of \eqref{eq:lex_example} is the one in \eqref{eq:na_example}.
Here, the objective function $\zeta(x)$ is built summing the first objective weighed by $\beta_1=1$ with the second one weighed by $\beta_2=\eta$, i.e., $\zeta(x)=f_1(x)+f_2(x)\eta$.
Notice that the request $\tfrac{\beta_2}{\beta_1}\approx 0$ holds since $\tfrac{\beta_2}{\beta_1}=\tfrac{\eta}{1}=\eta\approx 0$.
Furthermore, it is clear why the origin optimizes also the second problem.
Indeed, as soon as $x_1\ge0$ or $x_2\ge0$ then $\zeta(x)$ assumes positive finite values, while if $x_3\ge0$ then $\zeta(x)$ assumes positive infinitesimal values.
On the other hand, in the origin it holds $\zeta(x)=0$, which is a value smaller than any other assumed by the objective function in the previous two cases.\\
\vspace{2mm}
\begin{minipage}{.5\linewidth}
\vspace{3mm}
\begin{mini}|s|
{x\in\mr^3}{x_1^2+x_2^2,\,x_3+x_2}{\label{eq:lex_example}}{\text{\hspace{-1cm}lex}}
\addConstraint{0\le x_1\le1}{}{}
\addConstraint{0\le x_2\le1}{}{}
\addConstraint{0\le x_3\le1}{}{}
\end{mini}
\end{minipage}
\begin{minipage}{.5\linewidth}
\vspace{3mm}
\begin{mini}|s|
{x\in\mr^3}{x_1^2+x_2^2+(x_3+x_2)\eta}{\label{eq:na_example}}{}
\addConstraint{0\le x_1\le1}{}{}
\addConstraint{0\le x_2\le1}{}{}
\addConstraint{0\le x_3\le1}{}{}
\end{mini}
\end{minipage}

\section{Non-Archimedean Interior Point Method} 
\label{sec:naipm}
This section aims to introduce a non-Archimedean extension of the predictor-corrector primal-dual infeasible IPM described in Section \ref{sec:ipm}, called by the authors non-Archimedean IPM (NA-IPM).
Such algorithm is able to solve non-Archimedean quadratic optimization problems and, as a corner case, standard QP ones as well.
Three immediate and concrete advantages in leveraging it are: i) lexicographic QP problems can be solved as they were unprioritized ones (thanks to Proposition \ref{prop:convertion}); ii) the problem of infeasibility and unboundedness disappears; and iii) more general QP optimization problems can be modelled and solved (Alpha-Theory helps in modelling problems that are difficult to model without it, sometimes it is even impossible). 

The algorithm pseudocode is presented in Algorithm \ref{alg:na_IPM}. 
To better understand it, the next three subsections are devoted to discuss some delicate aspects of it and of its implementation: unboundedness and infeasibility, complexity and convergence properties, numerical issues which rise when moving from algorithmic description to software implementation.
Similarly to what happens for standard algorithms, it is important to stress that many of the proofs in this section assume Euclidean numbers to be represented by finite sequences of monosemia.
Indeed, even if the reference set $\me$ defined by the Axioms 1-3 admits numbers represented by infinite sequences, it would not be reasonable to use them in a machine and to discuss about algorithms convergence.
The reasons are two: i) the algorithm should manage and manipulate an infinite amount of data; ii) the machine is finite and cannot store all that information.
Notice that, at this stage, the focus is not on variable-length representations of Euclidean numbers, as they would slow the computations down \cite{Fiaschi_J4}.
Fixed-length representations, as the ones discussed in \cite{Benci_J1} are therefore preferred, also because they are easier to implement in hardware (i.e., they are more ``hardware friendly"), as recent studies testify \cite{bpu}.


\begin{algorithm}[!htbp]
\caption{Non-Archimedean predictor-corrector infeasible primal-dual IPM}
\label{alg:na_IPM}
\begin{algorithmic}[1]
\Procedure{NA-IPM}{$A$, $b$, $c$, $Q$, $\varepsilon$, $\texttt{max\_it}$}
\LineComment{Notice that the divergence is dealt with the embedding presented in Section \ref{sec:naipm_unboundedness}}
\LineComment{Therefore the flag of correct termination and the threshold $\omega$ are useless here}
\LineComment{Notice also that only $\varepsilon\in\mr^+$, while $A,b,c$ and $Q$ are Euclidean matrices and vectors}
\State $x$, $\lambda$, $s$ = starting\_point($A$, $b$, $c$, $Q$)
\State $x=\mathtt{lead\_mon}(x) \quad \lambda=\mathtt{lead\_mon}(\lambda) \quad s=\mathtt{lead\_mon}(s)$
\For{$i=1,\,\ldots,\,\texttt{max\_it}$}
    \LineComment{compute residuals}
    \State $r_b = Ax-b \quad r_c = A^T\lambda+s-Qx-c \quad r_{\mu} = x \odot s$
    \LineComment{compute centrality, $n$ is the length of $x$}
    \State $\mu = \tfrac{r_{\mu}^T\mathbf{1}}{n}$
    \LineComment{compute KKT conditions satisfaction parameters}
    \State $\rho_1 = \tfrac{\|Ax-b\|}{1\mathcal{O}(b)+\|b\|}, \quad
    \rho_2 = \tfrac{\|A^T\lambda+s-Qx-c\|}{1\mathcal{O}(c)+\|c\|}, \quad
    \rho_3 = \tfrac{\mu}{1\mathcal{O}(\frac{1}{2}x^TQx+c^Tx)+|\frac{1}{2}x^TQx+c^Tx|}$
    \LineComment{check convergence on all the (meaningful) monosemia of the aggregated values (Section \ref{sec:naipm_numerical})}
    \If{$\mathtt{all\_monosemia}(\rho_1)\le\varepsilon$ \text{and} $\mathtt{all\_monosemia}(\rho_2)\le\varepsilon$ \text{and} $\mathtt{all\_monosemia}(\rho_3)\le\varepsilon$}
        \LineComment{primal-dual feasible optimal solution found}
        \State \Return $x$, $\lambda$, $s$
    \EndIf
    \LineComment{compute predictor directions solving \eqref{eq:KKT}}
    \State $\Delta x_p$, $\Delta \lambda_p$, $\Delta s_p = $ predict($A$, $b$, $c$, $Q$, $r_b$, $r_c$, $r_{\mu}$)
    \LineComment{keep only the leading monosemia of the gradients (Section \ref{sec:naipm_numerical})}
    \State $\Delta x_p = \mathtt{lead\_mon}(\Delta x_p) \quad \Delta \lambda_p = \mathtt{lead\_mon}(\Delta \lambda_p) \quad \Delta s_p = \mathtt{lead\_mon}(\Delta s_p)$
    \LineComment{compute predictor step size}
    \State $\mhl{\nu_{pp} = 0.99\min(\max\nolimits_{\overline{\nu}}\{\mathtt{lead\_mon}(\overline{\nu})\,|\, x+\overline{\nu}\Delta x_p \ge 0\},\,1)}$
    \State $\mhl{\nu_{pd} = 0.99\min(\max\nolimits_{\overline{\nu}}\{\mathtt{lead\_mon}(\overline{\nu})\,|\, s+\overline{\nu}\Delta s_p \ge 0\},\,1)}$
    \State $\nu=\min(\nu_{pp},\,\nu_{pd})$
    \LineComment{estimate $\sigma$}
    \State $\Tilde{x} = x+\nu\Delta x_p \quad \Tilde{s} = s+\nu\Delta s_p$
    \State $\mu^{\text{new}} = \tfrac{\Tilde{x}^T\Tilde{s}}{n}$
    \State $\sigma = \mathtt{lead\_mon}((\tfrac{\mu^{\text{new}}}{\mu})^3$)
    \LineComment{compute corrector directions solving the corresponding Newton's system}
    \State $\Delta x_c$, $\Delta \lambda_c$, $\Delta s_c = $ corrector($A$, $b$, $c$, $Q$, $\sigma\mu\mathbf{1}-\Delta x_p \odot \Delta s_p$)
    \LineComment{compute new direction}
    \State $\Delta x = \Delta x_p + \Delta x_c \quad \Delta \lambda = \Delta \lambda_p + \Delta \lambda_c \quad \Delta s = \Delta s_p + \Delta s_c$
    \LineComment{keep only the leading monosemia of the gradients (Section \ref{sec:naipm_numerical})}
    \State $\Delta x = \mathtt{lead\_mon}(\Delta x) \quad \Delta \lambda = \mathtt{lead\_mon}(\Delta \lambda) \quad \Delta s = \mathtt{lead\_mon}(\Delta s)$
    \LineComment{compute step size}
    \State $\mhl{\nu_{p} = 0.99\min(\max\nolimits_{\overline{\nu}}\{\mathtt{lead\_mon}(\overline{\nu})\,|\, x+\overline{\nu}\Delta x \ge 0\},\,1)}$
    \State $\mhl{\nu_{d} = 0.99\min(\max\nolimits_{\overline{\nu}}\{\mathtt{lead\_mon}(\overline{\nu})\,|\, s+\overline{\nu}\Delta s \ge 0\},\,1)}$
    \State $\nu=\min(\nu_p,\,\nu_d)$
    \LineComment{compute target primal-dual solution}
    \State $x = x+\nu\Delta x \quad \lambda = \lambda+\nu\Delta \lambda \quad s = s+\nu\Delta s$
    \LineComment{Add infinitesimal centrality to the close-to-zero entries (Section \ref{sec:naipm_numerical})}
    \State $x,\,s = \mathtt{update\_zero\_entries}(x,\,s)$
\EndFor
\State \Return $x$, $\lambda$, $s$
\EndProcedure
\end{algorithmic}
\end{algorithm}

\subsection{Infeasibility and Unboundedness}
\label{sec:naipm_unboundedness}
As stated in Section \ref{sec:ipm}, one can approach the problem of infeasibility and unboundedness in two different ways: divergence detection at run-time or problem embedding.
While the first keeps the problem complexity fixed but negatively affects the computation because of norm divergence polling, the second wastes resource by optimizing a more complex problem, which is solved efficiently nevertheless.
Therefore, the simpler the embedding is, the lesser it affects the performance.

One very simple embedding, proposed by \cite{megiddo_IPM,kojima_IPM,monteiro_IPM_p1,monteiro_IPM_p2}, consists in the following mapping: 
\begin{equation}
\begin{split}
    A\mapsto\Tilde{A}=
    \begin{bmatrix}
    A & b-A\ones & \zeros\\
    c^T-\ones^T + \ones^TQ & 0 & -1
    \end{bmatrix}, \quad \quad \quad \quad &
    Q \mapsto \Tilde{Q}=
    \begin{bmatrix}
    Q & \zeros\\
    \zeros & \zeros
    \end{bmatrix},\\
    b\mapsto\Tilde{b}=
    \begin{bmatrix}
    b^T & -\wp_1
    \end{bmatrix}^T, \quad \quad \quad \quad &
    c\mapsto\Tilde{c}=
    \begin{bmatrix}
    c^T & \wp_2 & 0
    \end{bmatrix}^T,
\end{split}
\label{eq:standard_enlargement}
\end{equation}
where $\wp_1$ and $\wp_2$ are two positive and sufficiently big constants.
This embedding adds two artificial variables (one to the primal and one to the dual problem) and one slack variable (to the primal).
The goal of adding the artificial variables is to guarantee the feasibility of their corresponding problem, while on their own dual this is equivalent to add one bounding hyperplane to prevent any divergence.
From duality theory indeed, if the primal problem is infeasible, then the dual is unbounded and vice versa.
Geometrically, the hyperplane slope is chosen considering a particular conical combination of the constraints and the constant term vector (actually, the one with all coefficients equal to 1).
If there is any polyhedron unboundedness, the conical combination outputs a diverging direction and generates an hyperplane orthogonal to it; otherwise, the addition of such constraints has no effect.

On the other hand, the constraint intercept depends on the penalizing weights $\wp_1$ and $\wp_2$, respectively for the primal and the dual hyperplane.
The larger the weight is, the farther is located the corresponding bound.
On the primal perspective instead,
$\wp_1$ and $\wp_2$ act as penalizing weights for the artificial variables, respectively of the dual and the primal problem.
The need of this penalization comes from the fact that, to make the optimization consistent, the algorithm must be driven towards feasible points of the original problem, if any.
By construction, the latter always have artificial variables equal to zero, which means one has to penalize them in the cost function as much as possible in order to force them to that value.
More formally, it can be proved that for sufficiently large values of $\wp_1$ and $\wp_2$: i) the enlarged problem is strictly feasible and bounded; ii) any solution for the larger problem is also optimal for the embedded one if and only if both the artificial variables are zero \cite{kojima_IPM}.

Unfortunately, this idea is unsustainable when moving from theory to practice, i.e., to implementation.
Indeed, a good estimate of the weights is difficult to determine a priori, and the computational performance is sensitive to their values \cite{mcshane_IPM}.
Trying to figure out a solution, Lustig \cite{lustig_IPM_feasibility} investigated the optimal directions generated by Newton's step equation when $\wp_1$ and $\wp_2$ are driven to $\infty$, proposing a weight-free algorithm based on these directions.
Later, Lustig et al. \cite{lustig_91} showed that directions coincide with those of an infeasible IPM, without solving the unboundedness issue actually.
When considering a set of numbers larger than $\mr$ as $\me$ however, an approach in the middle between \eqref{eq:standard_enlargement} and the one by Lustig is possible.
It consists in the use of \emph{infinitely} large penalizing weights, i.e., in a non-Archimedean embedding.
This choice has the effect to infinitely penalize the artificial variables, while from a dual perspective it locates the bounding hyperplanes infinitely far from the origin.
For instance, in the case of a standard QP problem it is enough to set both $\wp_1$ and $\wp_2$ to $\alpha$, obtaining the following map
\begin{equation*}
    b\mapsto\Tilde{b}=
    \begin{bmatrix}
    b^T & -\alpha
    \end{bmatrix}^T, \quad \quad \quad \quad
    c\mapsto\Tilde{c}=
    \begin{bmatrix}
    c^T & \alpha & 0
    \end{bmatrix}^T.
\end{equation*}
The idea to infinitely penalize an artificial variable is not completely new: actually, it has already been successfully used in the I-Big-M method \cite{IBigM_20}, previously proposed by the author of this work, even if in a discrete context rather than in a continuous one.

Nevertheless, there is still a little detail to care of.
Embedding-based approaches leverage on the milestone theorem of duality to guarantee optimal solutions existence and boundedness.
A non-Archimedean version of the duality theorem must hold too, otherwise non-Archimedean embeddings end up to be theoretically not well-founded.
Thanks to the Transfer Principle, $\me$ is free from any issue of this kind, as stated by the next proposition-
\begin{myproposition}[Non-Archimedean Duality]
\label{theo:duality}
Given an NA-QP maximization problem, suppose that the primal and dual problems are feasible. Then, if the dual problem has a strictly feasible point, the optimal primal solution set is nonempty and bounded. The vice versa is true as well.
\begin{proof}
The theorem is true thanks to the Transfer Principle which, roughly speaking, transfers the properties of standard quadratic functions to quadratic non-Archimedean ones.
\end{proof}
\end{myproposition}
\noindent If a generic non-Archimedean QP problem is considered instead, setting the weights to $\alpha$ may be insufficient to correctly build the embedding.
Actually, their proper choice depends on the magnitude of the values constituting the problem.
Proposition \ref{prop:weights} gives a sufficient estimate of them; before showing it however, three preliminary results are necessary.
Lemma \ref{lem:bounds_on_optimals}, Lemma \ref{lem:bounds_unbounded} and Lemma \ref{lem:bounds_infeasible} address them. All these three lemmata make use of the functions $\mathcal{O}(\cdot)$ and $\mathit{o(\cdot)}$ provided in Section \ref{sec:at} as Definition 4 and 5, respectively.
In particular, Lemma \ref{lem:bounds_on_optimals} provides an upper bound to the magnitude of the entries of the solutions  $x$ of a non-Archimedean linear system $Ax\le b$.
This upper bound is expressed as function of the magnitude of the entries of both $A$ and $b$.
Furthermore, Lemma \ref{lem:bounds_on_optimals} considers the case in which the linear system to solve is the dual feasibility constraint of a QP problem, i.e., it has the form $A^T\lambda-Q\overline{x}\le c$ with $\overline{x}$ satisfying $A\overline{x}\le b$.
Lemmata \ref{lem:bounds_unbounded} and \ref{lem:bounds_infeasible} generalize  Lemma \ref{lem:bounds_on_optimals} considering corner cases too.
\begin{mylemma}
\label{lem:bounds_on_optimals}
Let the set of primal-dual optimal solutions $\Omega$ be nonempty and bounded.
Also, let $b,\,[c,\,Q]\neq\zeros$, $A$ has full row rank, and its entries are represented by at most $l$ monosemia, i.e., $A_{ij}=\sum_{h=1}^l (a_{ij})^h\alpha^{g(h)}$.
Then, any $(\overline{x},\,\overline{\lambda},\,\overline{s})\in\Omega$ satisfies 
\begin{equation*}
    \mathcal{O}(\overline{x}) \le \mathcal{O}(\Tilde{x}) = \min_{j\in J}\frac{\mathcal{O}(b_j)}{\mathit{o}(A_j)},
    \quad\quad
    \mathcal{O}(\overline{\lambda}) \le \mathcal{O}(\Tilde{\lambda}) = \min_{i\in I}\frac{\mathcal{O}([c_i,\,Q_i\mathcal{O}(\Tilde{x})])}{\mathit{o}([A^i,\, 1])},
\end{equation*}
where $J=\{j=1,\,\ldots,\,m\,|\,b_j\neq0\}$ and $I=\{i=1,\,\ldots,\,n\,|\,c_i\neq0\;\lor\;Q_i\neq\zeros\}$.
\begin{proof}
By hypothesis, $A\overline{x}=b$ and therefore $|A\overline{x}|=|b|$ too.
Focusing on the $j$-th constraint, $j\in J$, it holds
\begin{equation}
    |b_j| = \left|\sum_{h=1}^l (A_j)^h\overline{x}\alpha^{g(h)}\right| \ge |(A_j)^l\overline{x}\alpha^{g(l)}|,
\label{eq:orthogonal_simplification}
\end{equation}
which implies
\begin{equation}
    \mathcal{O}(b_j) \ge \mathcal{O}(|(A_j)^l\overline{x}\alpha^{g(l)}|) = \mathit{o}(A_j)\mathcal{O}(\overline{x}) \Longrightarrow \mathcal{O}(\overline{x}) \le \frac{\mathcal{O}(b_j)}{\mathit{o}(A_j)}\Longrightarrow \mathcal{O}(\overline{x}) \le \min_{j\in J} \frac{\mathcal{O}(b_j)}{\mathit{o}(A_j)}.
\label{eq:x_magn_upper_bound}
\end{equation}
The proof for the second part of the thesis is very similar:
\begin{equation*}
    A^T\overline{\lambda} + I\overline{s} = c+Q \overline{x} \Longrightarrow
    |A^{iT}\overline{\lambda} + \overline{s}_i| = |c_i +Q_i \overline {x}| \le |c_i|+|Q_i \overline{x}| \le |c_i|+|Q_i\xi|,
\end{equation*}
where $\xi\in\me^n$ is such that $|Q_i\overline{x}|\le|Q_i\xi|$ and has the form $\xi=\xi^0 \mathcal{O}(\Tilde{x})$, $\xi^0\in\mr^n$.
Now, following the same guidelines used in \eqref{eq:orthogonal_simplification} and \eqref{eq:x_magn_upper_bound}, one gets
\begin{equation*}
    \mathcal{O}(\overline{\lambda})\le
    \mathcal{O}([\overline{\lambda},\,\overline{s}]) \le
    \frac{\mathcal{O}([c_i,\,Q_i\mathcal{O}(\Tilde{x})])}{\mathit{o}([A^i,\, 1])} \quad \forall\,i\in I \Longrightarrow
    \mathcal{O}(\overline{\lambda})\le \min_{i\in I}\frac{\mathcal{O}([c_i,\,Q_i\mathcal{O}(\Tilde{x})])}{\mathit{o}([A^i,\, 1])}.
\end{equation*}
\end{proof}
\end{mylemma}
\noindent Equation \eqref{eq:orthogonal_simplification} may seem analytically trivial, but actually it underlines a subtle property of non-Archimedean linear systems: the solution can have entries infinitely larger than any number involved in the system itself.
As an example, the 2-by-2 linear system below admits the unique solution $\overline{x}=[\alpha^2,\,\alpha^2]$.
However, each value in the system is finite, i.e., the magnitude of each entry of $A$ and $b$ is $\mathcal{O}(\alpha^0)=\mathcal{O}(1)$:
\begin{equation*}
    \begin{bmatrix}
    \eta^2-1 & 1\\
    1 & \eta^2-1
    \end{bmatrix}
    x = 
    \begin{bmatrix}
    1\\1
    \end{bmatrix}
\end{equation*}
Notice that Lemma \ref{lem:bounds_on_optimals} works perfectly here.
Indeed, $\mathcal{O}(b)=1$ and $\mathit{o}(A)=\eta^2$ imply $\mathcal{O}(\Tilde{x})=\alpha^2\ge\mathcal{O}(\overline{x})$.

\begin{mylemma}
\label{lem:bounds_unbounded}
Let either the primal problem be unbounded or $\Omega\neq\emptyset$ be unbounded in the primal variable.
Let also $b,\,c,\,Q$ and $A$ satisfy the same hypothesis as in Lemma \ref{lem:bounds_on_optimals}.
Then,
\begin{equation*}
    \mathcal{O}(\overline{x}) \le \mathcal{O}(\Tilde{x}) = \min_{j\in J}\frac{\mathcal{O}(b_j)}{\mathit{o}(A_j)},
    \quad\quad
    \mathcal{O}(\overline{\lambda}) \le \mathcal{O}(\Tilde{\lambda}) = \min_{i\in I}\frac{\mathcal{O}([c_i,\,Q_i\mathcal{O}(\Tilde{x})])}{\mathit{o}([A^i,\,(c_i-1+Q_i^T\ones),\, 1])},
\end{equation*}
with $I$ and $J$ as in Lemma \ref{lem:bounds_on_optimals}.
\begin{proof}
If the primal problem is unbounded, it means that $\exists x\in\me^n$ such that $Ax=b$, and $\forall x$ such that $Ax=b,\,\nexists \lambda\in\me^m$ such that $A^T\lambda\le c+Q x$.
Nevertheless, a relaxed version of Lemma \ref{lem:bounds_on_optimals} still holds for the primal polyhedron, that is $\exists \overline{x} $ such that $A\overline{x}=b$ and $\mathcal{O}(\overline{x}) \le \mathcal{O}(\Tilde{x}) = \min_{j\in J}\tfrac{\mathcal{O}(b_j)}{\mathit{o}(A_j)}$ (the request for optimality is missing).
According to \eqref{eq:standard_enlargement}, a feasible bound for the primal polyhedron is $(c^T-\ones^T+\ones^T Q)x-\zeta=-\wp_2$, $\zeta\ge0$, provided a suitable choice of $\wp_2$.
Indeed, it can happen that a wrong value for $\wp_2$ turns the unbounded problem into an infeasible one.
This aspect shall be discussed in Proposition \ref{prop:weights}, which specifies $\wp_2$ as a function of $\overline{x}$.

Choice of $\wp_2$ apart, the addition of the bound to the primal polyhedron guarantees that $\exists x'=(x,\,\zeta)\in\me^{n+1}$ such that it is feasible for the bounded primal problem and $\exists\lambda'\in\me^{m+1}$ such that $A^T\lambda+(c^T-\ones^T+\ones^T Q)^T\xi\le c+Q x$ (remember that $\xi$ is the dual variable associated to the new constraint of the primal problem and that $\Tilde{Q}x'=Q x$).
Following the same reasoning used in Lemma \ref{lem:bounds_on_optimals}, one gets the second part of the thesis
\begin{equation*}
    \mathcal{O}(\overline{\lambda}) \le \mathcal{O}([\overline{\lambda}',\,\overline{s}]) \le \min_{i\in I}\frac{\mathcal{O}([c_i,\,Q_i\mathcal{O}(\Tilde{x})])}{\mathit{o}([A^i,\,(c_i-1+Q_i^T\ones),\, 1])}.
\end{equation*}

The case in which $\Omega\neq\emptyset$ and unbounded in the primal variable is very similar.
Together with the assumption $b,\,[c,\,Q]\neq\zeros$, it means that there is a plenty of (not strictly) feasible primal-dual optimal solutions, but there does not exists any with maximum centrality.
This fact negatively affects IPMs, since they move towards maximum centrality solutions.
Therefore, an IPM which tries to optimize such a problem will never converge to any point, even if there is a lot of optimal candidates.
To avoid this phenomenon, it is enough to bound such set of solutions with the addition of a further constraint to the primal problem (which has also the effect to guarantee the existence of the strictly feasible solutions missing in the dual polyhedron).
As a result, the very same considerations applied for problem unboundedness work in this case as well, leading exactly to the result in the thesis of this lemma.
\end{proof}
\end{mylemma}

\begin{mylemma}
\label{lem:bounds_infeasible}
Let either the primal problem be infeasible or $\Omega\neq\emptyset$ be unbounded in the dual variable.
Let also $b,\,c,\,Q$ and $A$ satisfy the same hypothesis as in Lemma \ref{lem:bounds_on_optimals}.
Then,
\begin{equation*}
    \mathcal{O}(\overline{x}) \le \mathcal{O}(\Tilde{x}) = \min_{j\in J}\frac{\mathcal{O}(b_j)}{\mathit{o}([A_j,\,b-A\ones])},
    \quad\quad
    \mathcal{O}(\overline{\lambda}) \le \mathcal{O}(\Tilde{\lambda}) = \min_{i\in I}\frac{\mathcal{O}([c_i,\,Q_i\mathcal{O}(\Tilde{x})])}{\mathit{o}([A^i,\, 1])},
\end{equation*}
with $I$ and $J$ as in Lemma \ref{lem:bounds_on_optimals}.
\begin{proof}
In this case $\nexists x$ such that $Ax=b$ but $\forall x\in\me^n$ $\exists\lambda\in\me^m$ such that $A^T\lambda\le c + Q x$.
Enlarging the primal problem in accordance to \eqref{eq:standard_enlargement}, one has that $\exists x'=(\overline{x},\,\zeta)\in\me^{n+1}$ ($\zeta\ge0$) such that $A\overline{x}+(b-A\ones)\zeta=b$.
In addition, it holds $\exists \overline{\lambda}$ such that $A^T\overline{\lambda}\le c+Q\overline{x}$, provided that $(b-A\ones)^T\overline{\lambda}\le \wp_1$ for some suitable choice of $\wp_1$ (which shall be discussed in Proposition \ref{prop:weights} as well).
By an analogous reasoning to the ones used in Lemma \ref{lem:bounds_on_optimals} and \ref{lem:bounds_unbounded}, the thesis immediately comes.

The case in which $\Omega\neq\emptyset$ and is unbounded on the dual variable works in the same but symmetric way of the complementary scenario discussed in Lemma \ref{lem:bounds_unbounded}.
Because of this, it implies the bounds stated in the thesis, while the proof is omitted for brevity.
\end{proof}
\end{mylemma}

\begin{myproposition}
\label{prop:weights}
Given an NA-QP problem and its embedding as defined in \eqref{eq:standard_enlargement}, a sufficient estimate of the penalizing weights is
\begin{equation*}
    \mathcal{O}(\wp_2) = \mathcal{O}\left(\alpha\min_{j\in J}\frac{\mathcal{O}(b_j)}{\mathit{o}([A_j,\,b-A\ones])}\right),
    \quad \quad
    \mathcal{O}(\wp_1) = \mathcal{O}\left(\alpha\min_{i\in I}\frac{\mathcal{O}([c_i,\,Q_i\mathcal{O}(\wp_2\eta)])}{\mathit{o}([A^i,\,(c_i-1+Q_i^T\ones),\, 1])}\right),
\end{equation*}
with $I$ and $J$ as in Lemma \ref{lem:bounds_on_optimals}.
In case $J=\emptyset$ then $\wp_2=\alpha$, while $I=\emptyset$ implies $\wp_1=\alpha$.
\begin{proof}
The extension to the quadratic case of Theorem 2.3 in \cite{kojima_IPM} (proof omitted for brevity) gives the following sufficient condition for $\wp_1$ and $\wp_2$, which holds true even in a non-Archimedean context thanks to the Transfer Principle:
\begin{equation}
    \wp_1 > (c^T-\ones^TA+\ones^TQ)^T\overline{x}, \quad \quad \wp_2 > (b-A\ones)^T\overline{\lambda},
    \label{eq:suff}
\end{equation}
where ($\overline{x}$, $\overline{\lambda}$, $\overline{s}$) is an optimal primal-dual solution of the original problem, if any.
A possible way to guarantee Equation \eqref{eq:suff} satisfaction is to choose $\wp_1$ and $\wp_2$ such that their magnitudes are infinitely higher than the right hand terms of the inequalities.
For instance one may set
\begin{equation*}
    \wp_1 \ge \alpha (c^T-\ones^TA+\ones^TQ)^T\overline{x}, \quad \quad \wp_2 \ge \alpha (b-A\ones)^T\overline{\lambda},
\end{equation*}
or more weakly
\begin{equation*}
    \mathcal{O}(\wp_1) \ge \mathcal{O}(\alpha (c^T-\ones^TA+\ones^TQ)^T\overline{x}), \quad \quad 
    \mathcal{O}(\wp_2) \ge \mathcal{O}(\alpha (b-A\ones)^T\overline{\lambda}).
\end{equation*}

In case $\Omega$ is nonempty and bounded, Lemma \ref{lem:bounds_on_optimals} holds and provides an estimate on the magnitude of both $\overline{x}$ and $\overline{\lambda}$.
In case either the primal problem is unbounded or $\Omega$ is unbounded in the primal variable, Lemma \ref{lem:bounds_unbounded} applies: the optimal solution is handcrafted by bounding the polyhedron, its magnitude is overestimated by $\Tilde{x}$, and \eqref{eq:suff} gives a clue for a feasible choice of $\wp_2$.
Similar considerations hold for the case of either primal problem infeasibility  or $\Omega$ unboundedness in the dual variable, where Lemma \ref{lem:bounds_infeasible} is used.

Corner cases are the scenarios where either $J=\emptyset$ or $I=\emptyset$.
Since $J=\emptyset$ implies $b=\zeros$, the primal problem is either unbounded or with unique feasible (and optimal) point $x=\zeros$.
In both cases, it is enough to set $\mathcal{O}(\Tilde{x})=1$.
Since $\zeros$ is a feasible solution, in case of unboundedness it must exist a feasible point with at least one finite entry and no infinite ones because of continuity.
In the other scenario, $\zeros$ is the optimal solution and therefore any finite vector is a suitable upper bound for it.
Analogous considerations hold for the case $I=\emptyset$, where a sufficient magnitude bound is $\mathcal{O}(\Tilde{\lambda})=1$.
\end{proof}
\end{myproposition}


\subsection{Convergence and complexity}

The main theoretical aspects to investigate in an iterative algorithm are convergence and complexity.
Notice that in case of non-Archimedean algorithms, the complexity of elementary operations (such as the sum) assumes their execution on non-Archimedean numbers, rather than on real ones.
Since theoretically NA-IPM is just an IPM able to work with numbers in $\me$, one first result on NA-IPM complexity comes straightforwardly thanks to Transfer Principle.
It is worth stressing that, as usual, Theorem \ref{theo:convergence} assumes to apply NA-IPM to a NA-QP problem whose optimal solutions set is non-empty and bounded.
\begin{mytheorem}[NA-IPM convergence]
\label{theo:convergence}
NA-IPM algorithm converges in $\mathcal{O}(n^2|\log\varepsilon|)$, where $n\in\mn$ is the primal space dimension and $\varepsilon\in\me^+$ is the optimality relative tolerance.
\begin{proof}
The theorem holds true because of the Transfer Principle.
\end{proof}
\end{mytheorem}

In spite of this result being remarkable, it is of no practical utility.
Indeed, the relative tolerance $\varepsilon$ may be not a finite value but an infinitesimal one, making the time needed to converge infinite.
However, under proper assumptions, also the finite time convergence can be guaranteed, as stated by Theorem \ref{theo:convergence_bis}.
Before showing it, some preliminary results are needed and are presented as lemmas.
In fact, Lemma \ref{lem:finite_decreasing} guarantees optimality improvement iteration by iteration, Lemma \ref{lem:upper_d_mu} provides a preliminary result used by Lemma \ref{lem:level_convergence} which proves the algorithm convergence on the leading monosemium.

\begin{mylemma}
\label{lem:finite_decreasing}
In NA-IPM, if $\sigma\in\mr^+$ then $\exists\, \Tilde{\nu}\in[0,\,1]\subset\mr$ such that $\mu^{(k+1)}\le(1-0.1\Tilde{\nu})\mu^{(k)}$ and $\|(r_b^{(k+1)},\,r_c^{(k+1)})\|\le(1-0.1\Tilde{\nu})\|(r_b^{(k)},\,r_c^{(k)})\|$.
\begin{proof}
Applying the Transfer Principle to Lemma 6.7 in \cite{wright_pd_IPM}, it holds true that
\begin{equation}
\Tilde{\nu} = \min\left(\frac{n\sigma}{C},\,\frac{\sigma(1-\gamma)}{C},\,\frac{0.49n}{C},\,1\right),
\label{eq:nu}
\end{equation}
where $C$ is a positive constant at most finite.
Equation \eqref{eq:nu} immediately implies $\Tilde{\nu}\in[0,\,1]\subset\me$ and $\mathcal{O}(\Tilde{\nu}) = \mathcal{O}(\sigma)$.
The assumption $\sigma\in\mr^+$ completes the proof.
\end{proof}
\end{mylemma}

\begin{mylemma}
\label{lem:upper_d_mu}
Let $d^{(k)}$ be the right hand term in \eqref{eq:Newton_IPM} at the $k$-th iteration, and $d_{\mu}$ the vector of its last $n$ entries.
If the temporary solution $(x^{(k)},\,\lambda^{(k)},\,s^{(k)})\in\mathcal{N}_{\um\infty}(\gamma,\,\beta)$ (see Lemma \ref{lem:level_convergence} for its definition), then $\|d_{\mu}\|<n\mu$.
\begin{proof}
By definition, $\|d_{\mu}\| = \sqrt{\sum_{i=1}^n(\sigma\mu-x_i s_i)^2}$.
Focusing on the radicand, one has
\begin{equation*}
\begin{split}
    \sum_{i=1}^n(\sigma\mu-x_i s_i)^2 &=
    \sum_{i=1}^n\sigma^2\mu^2-2\sigma\mu x_i s_i+(x_i s_i)^2 =
    n\sigma^2\mu^2-2\sigma\mu\sum_{i=1}^n x_i s_i+\sum_{i=1}^n(x_i s_i)^2 =\\
    &= n\sigma^2\mu^2-2n\sigma\mu^2+\sum_{i=1}^n(x_i s_i)^2 <
    n\sigma^2\mu^2-2n\sigma\mu^2+n^2\mu^2,
\end{split}
\end{equation*}
where the strict inequality comes from the fact that $x_i s_i\ge\gamma\mu$ by hypothesis, which implies $x_i s_i\le n\mu-\gamma\mu(n-1)< n\mu$.
Considering again the square root, the result comes straightforwardly:
\begin{equation*}
    \|d_{\mu}\| < \sqrt{n\sigma^2\mu^2-2n\sigma\mu^2+n^2\mu^2} = \mu \sqrt{n\sigma^2-2n\sigma+n^2} \le n\mu.
\end{equation*}
\end{proof}
\end{mylemma}

\begin{mylemma}
\label{lem:level_convergence}
Let $(x^{(0)},\,\lambda^{(0)},\,s^{(0)})$ be NA-IPM starting point, and $M^{(0)}\Delta^{(0)}=d^{(0)}$ be the compact form for Newton's step equation \eqref{eq:Newton_IPM} at the beginning of the optimization.
Let one rewrite the right hand term $d^{(k)}=(d^{(k)})^0+(d^{(k)})^1$, where $(d_i^{(k)})^0 = \mathtt{lead\_mon}(d_i^{(k)}).$
Call also $d_r$ the first $n+m$ entries of $(d_i^{(k)})^0$ and $d_{\mu}$ the last $n$, i.e., $(d_i^{(k)})^0=(d_r,\,d_{\mu})$. 
Then, $\exists k\in\mn$ such that $|d_{r_i}|\le\varepsilon\mathcal{O}(\|(r_b^{(0)},\,r_c^{(0)})\|)$ and $|d_{\mu_i}|< n\varepsilon\mathcal{O}(\|(r_b^{(0)},\,r_c^{(0)})\|)$ and the Newton-Raphson method reaches that iteration in $\mathcal{O}(n^2|\log\varepsilon|)$.
\begin{proof}
As usual, the central path neighborhood is
\begin{equation*}
\begin{split}
\mathcal{N}_{\um\infty}(\gamma,\,\beta) = \Big\{(x,\,\lambda,\,s)\in\me^{2n+m}\,&\Big|\,(x,\,s)>0,\, x_i s_i\ge\gamma\mu,\, \\ & \|(r_b^{(k)},\,r_c^{(k)})\|\le\frac{\|(r_b^{(0)},\,r_c^{(0)})\|}{\mu^{(0)}}\beta\mu^{(k)},\,\beta\ge1,\, \beta,\gamma\in\mr^+\Big\}.
\end{split}
\end{equation*}
Lemma \ref{lem:finite_decreasing} and $\mu$'s positivity implies that $\mathcal{O}(\mu^{(k+1)})\le\mathcal{O}(\mu^{(k)})$ $\forall k \in\mn$.
The application of the Transfer Principle to Theorem 6.2 in \cite{wright_pd_IPM} guarantees that $\exists k'\in\mn$ such that $\mu^{(k')}\le\varepsilon\mathcal{O}(\mu^{(0)})$ holds true and the Newton-Raphson algorithm reaches that iteration in $\mathcal{O}(n^2|\log\varepsilon|)$.
Together, Lemma \ref{lem:finite_decreasing} and $(x^{(k)},\,\lambda^{(k)},\,s^{(k)})\in\mathcal{N}_{\um\infty}(\gamma,\,\beta)$ guarantee that $\exists k''$ (reached in polynomial time as well) such that $\|(r_b^{(k'')},\,r_c^{(k'')})\|\le\varepsilon\mathcal{O}(\|(r_b^{(0)},\,r_c^{(0)})\|)$ too.
Set $k=\max(k',\,k'')$.
Then, one has $\mu^{(k)}\le\varepsilon\mathcal{O}(\mu^{(0)})\Longrightarrow \mathtt{lead\_mon}(\mu^{(k)})\le\varepsilon\mathcal{O}(\mu^{(0)})$ and $\|(r_b^{(k)},\,r_c^{(k)})\|\le\varepsilon\mathcal{O}(\|(r_b^{(0)},\,r_c^{(0)})\|) \Longrightarrow \mathtt{lead\_mon}(\|(r_b^{(k)},\,r_c^{(k)})\|)\le\varepsilon\mathcal{O}(\|(r_b^{(0)},\,r_c^{(0)})\|)$.
Moreover, by construction 
it holds 
$\|d_r^{(k)}\|\le\mathtt{lead\_mon}(\|(r_b^{(k)},\,r_c^{(k)})\|)$ and $\|d_{\mu}^{(k)}\|< n\mathtt{lead\_mon}(\mu^{(k)})$ (the latter comes from Lemma \ref{lem:upper_d_mu}).
Therefore the following two chains of inequalities hold true:
\begin{equation*}
\begin{gathered}
    |d_{r_i}^{(k)}|\le\|d_r^{(k)}\|\le\mathtt{lead\_mon}(\|(r_b^{(k)},\,r_c^{(k)})\|)\le\varepsilon\mathcal{O}(\|(r_b^{(0)},\,r_c^{(0)})\|),\\
    |d_{\mu_i}^{(k)}|\le\|d_{\mu}^{(k)}\|<n\mathtt{lead\_mon}(\mu^{(k)})\le n\varepsilon\mathcal{O}(\mu^{(0)}),
\end{gathered}    
\end{equation*}
as stated in the thesis.
\end{proof}
\end{mylemma}

\begin{mycorollary}
\label{cor:x_s_small}
Let $k$ satisfy Lemma \ref{lem:level_convergence}, then either
\begin{equation*}
    x_i < \sqrt{n\varepsilon}\mathcal{O}(x_i) \quad \lor \quad s_i < \sqrt{n\varepsilon}\mathcal{O}(s_i)\;\; \forall\,i=1,\,\ldots,\,n.
\end{equation*}
\begin{proof}
The result comes straightforwardly from three facts: i) $x_i s_i < n\mu$ $\forall i=1,\,\ldots,\,n$; ii) $\mu\le\varepsilon\mathcal{O}(\mu^{(0)})$; iii) the leading term of entry of $x$, $s$ and $\lambda$ is never zeroed since the full optimizing step is never taken (see lines 19-20 and 31-32 in Algorithm \ref{alg:std_IPM}).
\end{proof}
\end{mycorollary}

\noindent We are now ready to provide the convergence theorem for the NA-IPM.
\begin{mytheorem}[NA-IPM convergence]
\label{theo:convergence_bis}
NA-IPM converges to the solution of a NA-QP problem in $\mathcal{O}(l n^2 |\log\varepsilon|)$, where $n\in\mn$ is the primal space dimension, $\varepsilon\in\mr^+$ is the relative tolerance, and $l\in\mn$ is the number of consecutive monosemia used in the problem optimization.
\begin{proof}
For the sake of simplicity, assume to represent all the Euclidean numbers in the NA-QP problem by means of the same function of powers $g\colon\mn\rightarrow\mq$.
From the approximation up to $l$ consecutive monosemia, one can rewrite $\|(r_b^{(k)},\,r_c^{(k)})\|=\sum_{i=1}^l r_i^{(k)}\alpha^{g(i)}$ and $\mu^{(k)}=\sum_{i=1}^l\mu_i^{(k)}\alpha^{g(i)}$.
Lemma \ref{lem:level_convergence} guarantees that $\exists\, k$ for which $(d^{(k)})^0$ is $\varepsilon$-satisfied.
Now, update the temporary solution substituting each entry of $x$ and $s$ which satisfies Corollary \ref{cor:x_s_small} with any feasible value one order of magnitude smaller, e.g., $x_i$ satisfying Corollary \ref{cor:x_s_small} is replaced with a positive value of the order $\alpha^{g(j+1)}$, where $j\in\mn$ is such that $\mathcal{O}(x_i)=\alpha^{g(j)}$.
Actually, they are those variables which are not active at the optimal solution, at least considering the zeroing of $(d^{(k)})^0$ only.
Then, recompute $d^{(k)}$, which by construction satisfies $\|d_{\mu}^{(k)}\|=0$, $\|d_r^{(k)}\|<\varepsilon\mathcal{O}(\|(r_b^{(0)},\,r_c^{(0)})\|)$.
All these operations have polynomial complexity and do not affect the overall result.
Updating the right hand term as $d^{(k)}\leftarrow d^{(k)}-(d^{(k)})^0$ and zeroing the leading term of those entries whose magnitude is still $\mathcal{O}(\|(r_b^{(0)},\,r_c^{(0)})\|)$, next algorithm iterations are forced to consider the previous $(d^{(k)})^0$ as already fully satisfied.
What is actually happening is that the problem now tolerates an infeasibility error whose norm is equal to $\varepsilon\mathcal{O}(\|(r_b^{(0)},\,r_c^{(0)})\|)$.
Therefore, one can apply Lemma \ref{lem:level_convergence} again to obtain one solution which is $\varepsilon$-optimal on the second monosemia of $(d^{(0)})^0$ too, and this result is achieved with a finite number of iterations and in polynomial complexity.
Repeating the update-optimization procedure for all the $l$ monosemia by means of which $(d^{(0)})^0$ is represented, one obtains one $\varepsilon$-optimal solution on all of them.
Since each of the $l$ $\varepsilon$-satisfactions is achieved in $\mathcal{O}(n^2|\log\varepsilon|)$, then the whole algorithm converges in $\mathcal{O}(l n^2|\log\varepsilon|)$.
\end{proof}
\end{mytheorem}

The next proposition highlights a particular property of NA-IPM when solving lexicographic QP problems.
Actually, it happens that every time $\mu$ decreases of one order of magnitude, then one objective is $\varepsilon$-optimized.
\begin{myproposition}
\label{prop:lex_convergence}
Consider an NA-QP problem generated from a standard lexicographic one in accordance with Theorem \ref{prop:convertion} and $\beta_i=\alpha^{1-i}$ $\forall i=1,\,\ldots,\,l$.
Then, each of the $l$ objectives is $\varepsilon$-optimized in polynomial time and when the $i$-th one is $\varepsilon$-optimized the magnitude of $\mu$ decreases from $\mathcal{O}(\alpha^{1-i})$ to $\mathcal{O}(\alpha^{-i})$ in the next iteration.
\begin{proof}
Assume to start the algorithm with a sufficiently good and well centered solution, as the one produced by Algorithm \ref{alg:start}, then, $\mathcal{O}(\|(r_b^{(0)},\,r_c^{(0)})\|)=\mathcal{O}(\mu^{(0)})=\mathcal{O}(\alpha^0)$.
Since $x^{(k)}\in\mr^+$ by construction, one can interpret each monosemia in $r_c^{(k)}$ as the satisfaction of the corresponding objective function at the $k$-th iteration, that is if $r_c^{(k)} = \sum_{i=1}^l r_{c_i}\alpha^{1-i}$ then the first objective lacks $r_{c_1}$ to be fully optimized, the second one lacks $r_{c_2}$ and so on.
Because of Lemma \ref{lem:level_convergence}, in polynomial time $(d_i)^0$ is $\varepsilon$-optimized, that is the KKT conditions \eqref{eq:KKT} are $\varepsilon$-satisfied.
In fact, this means that primal-dual feasibility is close to finite satisfaction and centrality is finitely close to zero.
There is a further interpretation nevertheless.
Interpreting the KKT conditions from a primal perspective their $\varepsilon$-satisfaction testifies that: i) the primal solution if feasible (indeed, the primal is a standard polyhedron and therefore is enough to consider the leading terms of $x$ only, getting rid of the infinitesimal infeasibility $r_{b_2}$); ii) the objective function is finitely $\varepsilon$-optimized (which means that the first objective is $\varepsilon$-optimized since the original problem was a lexicographic one and the high-priority objective is the only one associated to finite values of the non-Archimedean objective function); iii) the approximated solution is very close to the optimal surface of the first objective, roughly speaking it is $\varepsilon$-finitely close.
Moreover, the fact that $\|d_{\mu}\|=0$ after the updating procedure used in in Theorem $\ref{theo:convergence_bis}$  implies that the magnitude of $\mu$ will be one order of magnitude smaller in the next iteration, i.e., it will decrease from $\mathcal{O}(\alpha^0)$ to $\mathcal{O}(\alpha^{\um1})$.
Since what just said holds for all the $l$ monosemia (read priority levels of, i.e., objectives in the lexicographic cost function), the proposition is proved true.
\end{proof}
\end{myproposition}

\subsection{Numerical considerations and implementation issues}
\label{sec:naipm_numerical}
The whole field $\me$ cannot be used in practice, since it is too big to fit in a machine.
However, the algorithmic field $\widehat{\me}$ presented below is enough to represent and solve many real world problems:
\begin{equation*}
\widehat{\me}=\left\{\xi\in\me\,\Bigg|\, \xi = \sum_{i=1}^l\psi_i\alpha^{g(i)},\, l\in\mn,\,\psi_i\in\mr\right\}\cup\left\{0\right\},
\end{equation*}
where $g:\mn\rightarrow\mz$ is a monotone decreasing function and the term ``algorithmic field'' refers to finite approximations of theoretical fields realized by computers \cite{Benci_J2}.
Similarly to IEEE-754 floating point numbers which is the standard encoding for real numbers within a machine, a finite dimension encoding for Euclidean numbers in $\widehat{\me}$ is needed.
In \cite{Benci_J1,Benci_J2} the Bounded Algorithmic Number (BAN) representation is presented as a sufficiently flexible and informative encoding to 
cope with this task. The BAN format is a fixed-length approximation of an Euclidean number.
An example of BAN is $\alpha^{-1}(2.4 + 3.9\eta-2.89\eta^2)$, where the ``precision" in this context is given by the degree of the polynomial in $\eta$ plus 1 (three in this case). The BAN encoding with degree three is indicated as BAN3.

The second detail to care of when attempting to do numerical computations with Euclidean numbers is the effect of lower magnitude monosemia on them.
For instance, consider a two-objective lexicographic QP whose first objective is degenerate with respect to some entries of $x$.
When solving the problem by means of NA-IPM, the following phenomenon (which can also be proved theoretically) occurs: the information of the optimizing direction for the secondary objective is stored as an infinitesimal gradient in the solution of Newton's step equation \eqref{eq:Newton_IPM}.
As an example, assume that $x\in\mr^3$ and the entries $x_2$ and $x_3$ are degenerate with respect to the first objective.
Then, at each iteration the infinitesimal monosemium in the optimizing direction of $x_1$ assumes a negligible value, while for $x_2$ and $x_3$ this is not true: it is significant and grows exponentially in time.
In fact, the infinitesimal gradient represents the optimizing direction which must be followed along the optimal (and degenerate) surface of the first objective in order to reach optimality also for the second one.
However, such infinitesimal directions do not significantly contribute to the optimization, since the major role is played by the finite entries of the gradient.
Therefore, the effect of this infinitesimal information in the gradient only generates numerical instabilities.
As soon as the first objective is $\varepsilon$-optimized nevertheless, i.e., the first objective surface is reached, the optimizing direction still assumes finite values but this time oriented in order to optimize the second objective keeping the first one fixed, while all the infinitesimal monosoemia of the gradient assume negligible values. 
Roughly speaking, it happens a sort of ``gradient promotion'' as a result of the change in the objective to optimize.
To cope with the issue of noisy and unstable infinitesimal entries in the gradient, two details need to be implemented: i) after the computation of the gradients (both the predictor and the corrector step), only the leading term of each entry must be preserved, zeroing the remaining monosemia; ii) after having computed the starting point according to Algorithm \ref{alg:start}, again only the leading term of each entry of $x$, $s$, $\lambda$ must be preserved.
These variations do not affect convergence nor the generality of the discussion since the leading terms of the primal-dual solution are the only ones which impact on the zeroing of $(d^{(k)})^0$.

The choice of dealing with only the leading terms of the gradients comes in handy to solve another issue during the computations: a good choice for the value to assign to the zeroed entries of $x$ and $s$ during the updating phase discussed in Theorem \ref{theo:convergence_bis}.
Actually, it is enough to add them one monosemium whose magnitude is such that the following equality holds true: $\mathcal{O}(x_i s_i)=\mathcal{O}(\mu)\eta$.
For instance, assume again a two-objective lexicographic QP scenario after having completed the optimization of the first one.
It holds true that $\mathcal{O}(\mu)=\alpha^0$ and either $x_i$ or $s_i$ are smaller than $n\varepsilon\alpha^0$, say $x_i$.
Then, the updating phase of Theorem \ref{theo:convergence_bis} sets $x_i$ to a value having magnitude one order smaller.
A reasonable approach is to set $x_i$ equal to a monosemium, say $\xi$, such that $\mathcal{O}(\xi s_i)=\mathcal{O}(\mu)\eta=\eta$.
Since $\mathcal{O}(s_i)$ is finite because of Corollary \ref{cor:x_s_small}, one has $\mathcal{O}(\xi)=\eta$.
The naive choice is $\xi=\eta$, but it may be not the best one.
Indeed, this approach does not guarantee to generate a temporary solution with the highest possible centrality, i.e., $x_i s_i = \mu'$, where $\mu'$ is the centrality measure after the update.
To do it, the value to opt for must be $\xi_i = \tfrac{\mu'}{z_i}$, where $z_i = \max(x_i,s_i)$ and $\mu'$ is a monosemium one order of magnitude smaller than $\mu$ and sufficiently (but arbitrarily) far from zero.

Another numerical issue to care about is the computation accuracy due to the number of monosemia used.
As a practical example, consider the task to invert the matrix $A$ reported below.
Since its entries have magnitudes from $\alpha$ to $\eta$, one may think BAN3 encoding is enough to properly compute an approximation of its inverse, which is reported below as $A_3^{\um1}$.
However, the product $AA_3^{\um1}$ testifies the presence of an error whose magnitude is $\eta^2$, quite close to $\mathit{o}(A)=\eta$.
Depending on the application, such noise could negatively affect further computations and the final result.
Therefore, a good practice is to foresee some additional slots in the non-Archimedean numbers encoding.
For instance, adopting the BAN5 standard, the approximated inverse $A_5^{\um1}$ manifests an error with magnitude $\eta^4$ (see matrix $AA_5^{\um1}$), definitely a safer choice even if at the expenses of an extra computational effort.
\begin{equation*}
    A =
    \begin{bmatrix}
    \alpha & -\alpha & 2\eta \\ 
    2\alpha & \eta & -\alpha \\ 
    \eta & 2\alpha & -\alpha 
    \end{bmatrix},
    \quad \quad
    A_3^{\um1} =
    \begin{bmatrix}
    0.25 \alpha - 0.125 \eta & -0.12 \alpha + 0.5 \eta & 0.12 \alpha \\
    0.25 \alpha - 0.125 \eta & -0.12 \alpha & 0.12 \alpha + 0.5 \eta \\
    0.50 \alpha & -0.25 \alpha - 0.125 \eta & 0.25 \alpha + 0.125 \eta
    \end{bmatrix},
\end{equation*}
\begin{equation*}
    A_5^{\um1} = 
    \begin{bmatrix}
    0.25 \alpha - 0.125 \eta + 0.0625 \eta^{3} & -0.12 \alpha + 0.5 \eta - 0.0312 \eta^{3} & 0.12 \alpha - 0.219 \eta^{3}  \\
0.25 \alpha - 0.125 \eta + 0.0625 \eta^{3} & -0.12 \alpha - 0.281 \eta^{3} & 0.12 \alpha + 0.5 \eta + 0.0312 \eta^{3}  \\
0.5 \alpha & -0.25 \alpha - 0.125 \eta - 0.0625 \eta^{3} & 0.25 \alpha + 0.125 \eta + 0.0625 \eta^{3}
\end{bmatrix}
\end{equation*}
\begin{equation*}
    AA_3^{\um1} =
    \begin{bmatrix}
    1 & -0.25 \eta^{2} & 0.25 \eta^{2} \\
    0 & 1 & 0 \\
    0 & 0 & 1
    \end{bmatrix},
    \quad \quad
    AA_5^{\um1} =
    \begin{bmatrix}
    1 & -0.12 \eta^{4} & 0.12 \eta^{4}  \\
    0 & 1 & 0  \\
    0 & 0 & 1 
\end{bmatrix}
\end{equation*}


Finally, the last detail concerns terminal conditions.
In standard IPM, execution stops when the three convergence measures $\rho_1$, $\rho_2$ and $\rho_3$ in \eqref{eq:rho_def} are smaller than the threshold $\varepsilon$.
However, in a non-Archimedean context optimality means to be $\varepsilon$-optimal on each of the $l$ monosemia in objective function, i.e., to satisfy the KKT conditions \eqref{eq:KKT} on the first $l$ monosemia in $b$ and $c$ and $\mu$.
This means that the terminal condition needs to be modified in order to cope with such a convergence definition.
The convergence measures need to be redefined as
\begin{equation}
\begin{gathered}
\rho_1 = \frac{\|Ax-b\|}{\mathcal{O}(b)+\|b\|}, \quad\quad
\rho_2 = \frac{\|A^T\lambda+s-Qx-c\|}{\mathcal{O}(c)+\|c\|}, \\
\rho_3 = \frac{\mu}{\mathcal{O}(\frac{1}{2}x^TQx+c^Tx)+|\frac{1}{2}x^TQx+c^Tx|},
\end{gathered}
\label{eq:na_rho_def}
\end{equation}
with the convention $\mathcal{O}(0)=\alpha^0=1$.
In this way, one has the guarantee that $\rho_1$, $\rho_2$ and $\rho_3$ are finite values when close to optimality.
To better clarify this concept, consider the case in which $b$ has only infinitesimal entries.
This implies that the norm of $b$ is infinitesimal too.
In case the convergence measures in \eqref{eq:rho_def} are used, the denominator of $\rho_1$ is a finite number.
Therefore, any finite approximation of $b$, i.e., any primal solution $x$ such that $\|Ax-b\|$ is finite induces a finite value for $\rho_1$.
This is definitely a bad behavior since it is natural to expect that: i) the convergence measures are finite numbers only if the optimization is close to optimality; and ii) their leading monosemium is smaller than $\varepsilon$ when the leading monosemium of the residual norm ($\|Ax-b\|$ in this case) is small as well.
However, this is not the case of the current example as $\rho_1$ assumes finite values for approximation errors of $b$ which are infinitely larger than its norm.
Using the definitions in \eqref{eq:na_rho_def} instead, the issue is solved since now finite approximations of $b$ are mapped into infinite values of $\rho_1$, while infinitesimal errors into finite ones.
In fact, the introduction of the magnitude of the constant terms vectors in the definitions avoids the bias which would have been introduced by an a-priori choice of the magnitude of the constant term added to the denominator of the convergence measures.

Then, feasibility on the primal problem is $\varepsilon$-achieved when the absolute value of the first $l_b$ monosemia in $\rho_1$ are smaller than $\varepsilon$, i.e., assuming $b = \sum_{i=1}^{l_b} b^i\alpha^{g(i)}$ and $\rho_1 = \sum_{i=1}^l \rho_1^i\alpha^{1-i}$, a sufficient level of optimality is reached when $|\rho_1^i|<\varepsilon$ $\forall i=1,\,\ldots,\,l_b$, $l\ge l_b$.
Similar considerations hold for $c$ and $\mu$ too.

\section{Numerical Experiments}
\label{sec:experiments}
This section presents four problems by means of which the efficacy of NA-IPM is tested.
The first, the third and the fourth are lexicographic ones, while the second involves non-Archimedean numbers also in the constant term $b$.
They are listed in increasing order of difficulty, from LP to QP problems, from two to three objective functions, passing through unboundedness tackling.
In all the experiments $\varepsilon$ is set to $10^{\um8}$, a de facto standard choice.

\subsection{Experiment 1: two-objective LP}
The first experiment uses a benchmark already exploited to evaluate the efficacy of non-Archimedean optimizers \cite{CococcioniEtAlAMC2018,IBigM_20}.
In \cite{CococcioniEtAlAMC2018}, for instance, the algorithm used was a non-Archimedean version of the Simplex method, able to deal with non-Archimedean cost functions and standard constraints.
On the contrary, here the adopted algorithm is a non-Archimedean generalization of the primal-dual IPM, which is intrinsically implemented to cope with non-Archimedean functions both at the level of cost function and constraints.
Another difference consists in the fact that Simplex-like algorithms are discrete ones, while IPMs are continuous procedures.

\begin{figure*}[ht]
\begin{minipage}{.48\linewidth}
    \centering
    \includegraphics[width=\linewidth]{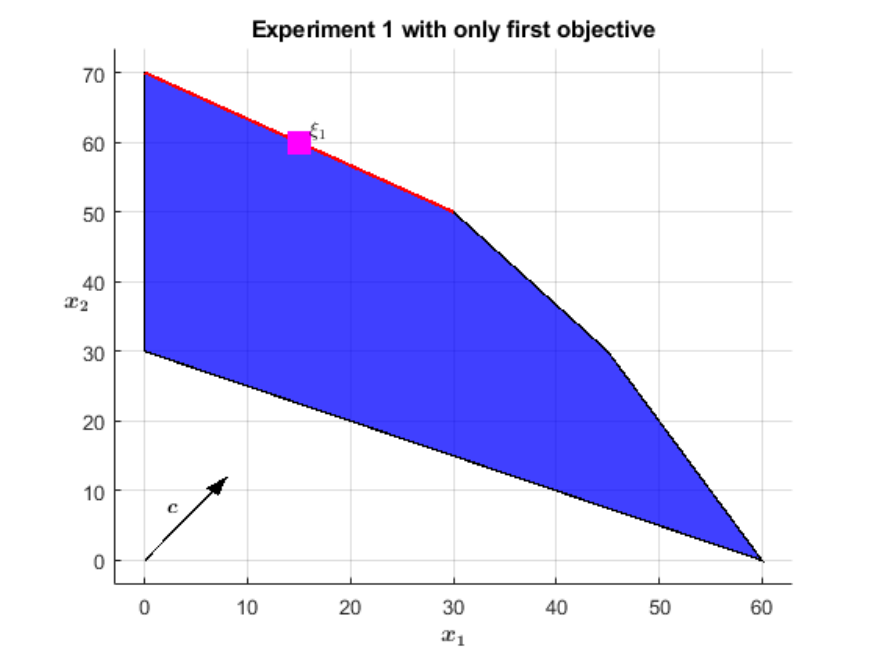}
    \caption{The optimal segment for the primary objective is in red; $\xi_1$ is the optimal point a standard IPM would approach.}
    \label{fig:ex_1_1}
\end{minipage}
\begin{minipage}{.04\linewidth}\end{minipage}
\begin{minipage}{.48\linewidth}
    \centering
    \includegraphics[width=\linewidth]{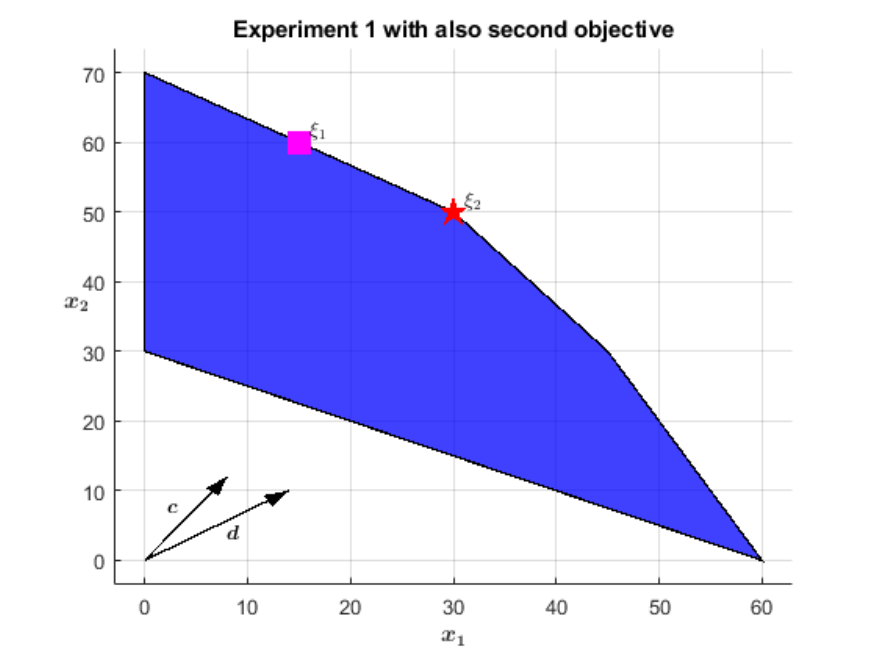}
    \caption{$\xi_2$ is the only optimal solution when both the objectives are considered.}
    \label{fig:ex_1_2}
\end{minipage}
\end{figure*}

The problem formulation is in Equation \eqref{eq:kyte}.
Geometrically, the first and most important objective $c = [8,\,12]^T$ identifies as optimal the whole segment linking the point $(0,\,70)$ to $(30,\,50)$, highlighted in red in Figure \ref{fig:ex_1_1}.
When considering also the second objective $d=[14,\,10]^T$, only one point of this segment turns to be truly optimal: the vertex $\xi_2 = (30,\,50)$.
A red star indicates it in Figure \ref{fig:ex_1_2}.
Equation \eqref{eq:na_kyte} reports the non-Archimedean representation of the problem; the two objectives are scalarized into one single cost function by a weighting sum: the first one is left as is while the second is scaled down by the factor $\eta$, testifying the precedence (read importance) relation between the two (recall from Section \ref{sec:at} that $\eta$ is an infinitesimal value, defined as the reciprocal of $\alpha$).\\
\vspace{2mm}
\begin{minipage}{.5\linewidth}
    \vspace{3mm}
    \begin{maxi}|s|
    {x}{8x_1+12x_2,\,14x_1+10x_2}{\label{eq:kyte}}{\text{\hspace{-1cm}lex}}
    \addConstraint{2x_1+\;\,x_2}{\le120}{}
    \addConstraint{2x_1+3x_2}{\le210}{}
    \addConstraint{4x_1+3x_2}{\le270}{}
    \addConstraint{x_1+2x_2}{\ge60}{}
    \addConstraint{x_1,\,x_2}{\ge0}{}
\end{maxi}
\end{minipage}
\begin{minipage}{.5\linewidth}
    \vspace{3mm}
    \begin{maxi}|s|
    {x}{(8+14\eta)x_1+(12+10\eta)x_2}{\label{eq:na_kyte}}{}
    \addConstraint{2x_1+\;\,x_2}{\le120}{}
    \addConstraint{2x_1+3x_2}{\le210}{}
    \addConstraint{4x_1+3x_2}{\le270}{}
    \addConstraint{x_1+2x_2}{\ge60}{}
    \addConstraint{x_1,\,x_2}{\ge0}{}
\end{maxi}
\end{minipage}

As already stressed, any IPM favors solutions with higher centrality and NA-IPM is not an exception.
This aspect becomes particularly evident in case of multiple optimal solutions, since the optimizing algorithm converges to the analytic centre of the optimal surface.
Proposition \ref{prop:lex_convergence} implicitly states that in lexicographic problems NA-IPM primarily optimizes the first objective, than the second one, and so on.
Since in \eqref{eq:kyte} the optimal region for the first objective is a segment,
in the first place NA-IPM moves towards its midpoint $\xi_1=(15,\,60)$, marked by a magenta square in Figure \ref{fig:ex_1_1}.
Eventually, the first cost function is $\varepsilon$-optimized and the secondary objective starts to mainly condition the optimization, moving the temporary solution towards $\xi_2$ until termination conditions are satisfied.
This change in leadership is testified when the centrality measure turns from finite to infinitesimal, as discussed in Theorem \ref{theo:convergence_bis}, Proposition \ref{prop:lex_convergence} and Section \ref{sec:naipm_numerical}.
With reference to Table \ref{tab:kyte_iterations}, one can appreciate this phenomenon between line 5 and 6.
At line 5, the current solution is very close to $\xi_1$, actually $(18.45,\,57.70)$, while $\mu$ has assumed for the last time a (very small) finite value.
These facts testify that the first objective is $\varepsilon$-satisfied, therefore in the next step the algorithm will move towards $\xi_2$ and show an infinitesimal centrality measure, as confirmed by line 6 where the target point is $(29.88,\,50.08)$.
Moreover, the improvement in the objective function from step 5 to step 6 is infinitesimal, which further stresses that the algorithm has already reached first objective $\varepsilon$-satisfaction and it is moving along its degenerate surface.
Similar numerical confirmations come from the other experiments too and are highlighted in boldface in the associated tables.

\begin{table}[!ht]
        \centering
        \caption{Iterations of NA-IPM solving the problem in \eqref{eq:kyte}}
        \begin{tabular}{|c|c|c|c|}
        \hline
        \textbf{iter} & $\bm{\mu} \in \mr$ & $\bm{x} \in \mr^2$ & $\bm{f(x)} \in \me$\\
        \hline
        0 & $273.00 $ & $\begin{bmatrix}98.80 & 40.51 \end{bmatrix}$ & $-1276.48  - 1.79e3 \eta$ \\
        \hline
        1 & $38.64 $ & $\begin{bmatrix}26.94 & 43.47 \end{bmatrix}$ & $-737.22  - 8.12e2 \eta$ \\
        \hline
        2 & $2.97 $ & $\begin{bmatrix}18.53 & 57.56 \end{bmatrix}$ & $-838.97  - 8.35e2 \eta$ \\
        \hline
        3 & $0.03 $ & $\begin{bmatrix}18.45 & 57.70 \end{bmatrix}$ & $-839.99  - 8.35e2 \eta$ \\
        \hline
        4 & $29.81e\um4 $ & $\begin{bmatrix}18.45 & 57.70 \end{bmatrix}$ & $-840.00  - 8.35e2 \eta$ \\
        \hline
        \textbf{5} & $\mathbf{2.82e\um6} $ & $\begin{bmatrix}\bm{18.45} & \bm{57.70} \end{bmatrix}$ & $\bm{-840.00  - 8.35e2 \eta}$ \\
        \hline
        \textbf{6} & $\bm{12.82 \eta}$ & $\begin{bmatrix}\mathbf{29.88} & \mathbf{50.08} \end{bmatrix}$ & $\bm{-840.00  - 9.19e2 \eta}$ \\
        \hline
        7 & $0.14 \eta$ & $\begin{bmatrix}30.00 & 50.00 \end{bmatrix}$ & $-840.00  - 9.20e2 \eta$ \\
        \hline
        8 & $1.40e\um3\eta$ & $\begin{bmatrix}30.00 & 50.00 \end{bmatrix}$ & $-840.00  - 9.20e2 \eta$ \\
        \hline
        9 & $1.41e\um5 \eta$ & $\begin{bmatrix}30.00 & 50.00 \end{bmatrix}$ & $-840.00  - 9.20e2 \eta$ \\
        \hline
        10 & $4.30e\um8 \eta$ & $\begin{bmatrix}30.00 & 50.00 \end{bmatrix}$ & $-840.00  - 9.20e2 \eta$ \\
        \hline
    \end{tabular}
\label{tab:kyte_iterations}
\end{table}

\subsection{Experiment 2: unbounded problem}
\label{sec:ex_2}
\begin{figure*}[ht]
\begin{minipage}{.5\linewidth}
    \centering
    \includegraphics[width=\linewidth]{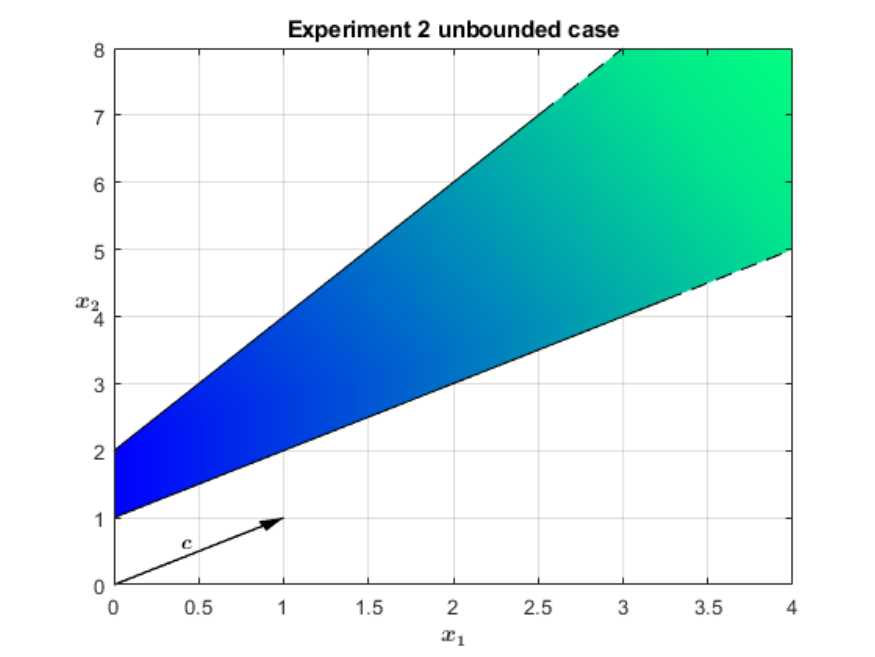}
    \caption{Example of unbounded primal polyhedron.}
    \label{fig:ex_2_1}
\end{minipage}
\begin{minipage}{.5\linewidth}
    \centering
    \includegraphics[width=\linewidth]{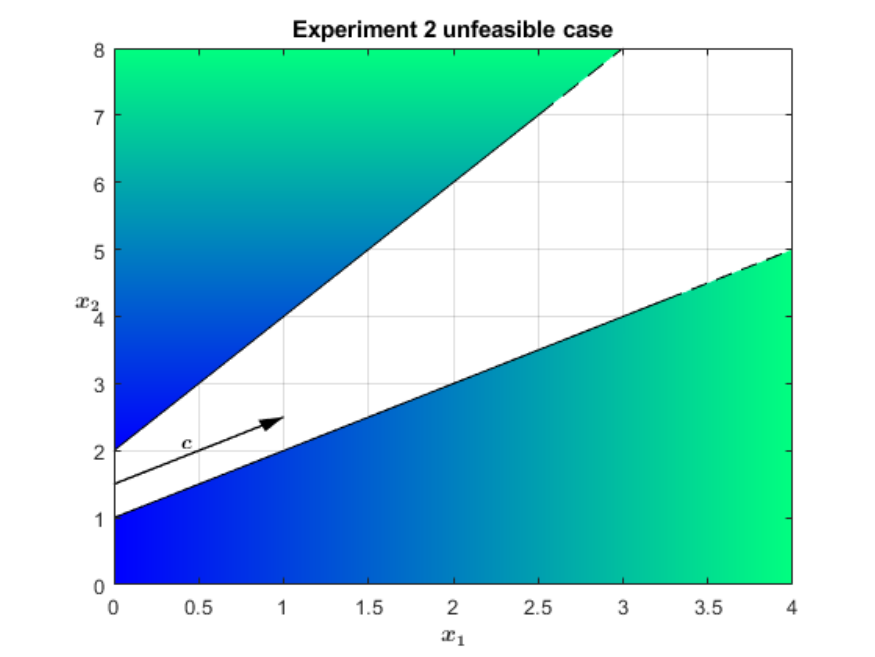}
    \caption{Example of empty primal polyhedron.}
    \label{fig:ex_2_2}
\end{minipage}
\end{figure*}

The second experiment aims to numerically show the efficacy of the mild embedding shown in Section \ref{sec:naipm_unboundedness} to cope with infeasibility and unboundedness.
As an example, consider the 2D unbounded problem described in Equation \eqref{eq:unbounded1} and drawn in Figure \ref{fig:ex_2_1}, which is already analytically reported in normal form as in \eqref{eq:QP_normal} for the sake of clarity.
To mitigate the issues coming from the iterates divergence, one can resort to the embedding described in Equation \eqref{eq:standard_enlargement}, obtaining the strictly feasible and bounded problem in \eqref{eq:unbounded1_embed}.
Proposition \ref{prop:weights} recommends the use of penalizing weights such that $\mathcal{O}(\wp_1)=\mathcal{O}(\wp_2)=\mathcal{O}(\alpha)$; the choice has been $\wp_1=\wp_2=\alpha$.\\
\vspace{2mm}
\begin{minipage}{.5\linewidth}
    \vspace{3mm}
    \begin{maxi}|s|
    {x}{\quad\; x_1+x_2}{\label{eq:unbounded1}}{}
    \addConstraint{-2x_1+x_2+x_3}{=2}{}
    \addConstraint{x_1-2x_2+x_4}{=1}{}
    \addConstraint{x}{\ge0}{}
    \addConstraint{x}{\in \mr^4}{}
    \end{maxi}
\end{minipage}
\begin{minipage}{.5\linewidth}
    \vspace{3mm}
    \begin{maxi}|s|
    {x}{\quad\; x_1+x_2-\alpha x_5}{\label{eq:unbounded1_embed}}{}
    \addConstraint{-2x_1+x_2+x_3+2x_5}{=2}{}
    \addConstraint{x_1-2x_2+x_4+x_5}{=1}{}
    \addConstraint{-x_3-x_4-x_6}{=-\alpha}{}
    \addConstraint{x}{\ge0}{}
    \addConstraint{x}{\in \me^6}{}
    \end{maxi}
\end{minipage}

Table \ref{tab:unbounded_iterations} reports the iterations made by NA-IPM to solve such an extended problem.
As expected, the algorithm converges in a finite number of steps, and the optimal point lies on the bounding hyperplane $-x_3-x_4-x_6=-x_1-x_2-3-x_6=-\alpha$ located infinitely far from the origin.
Formally, what gives clue about the unboundedness of the problem is the dual variable $\lambda_3$, see Proposition \ref{prop:weights}.
If the problem is bounded then it must be zero in the optimal solution, while it is equal to $1$.
In this specific case however, there is another and more significant indicator: the magnitude of $x_1$ and $x_2$.
Since the problem was a standard one before the embedding, if its solution exists it must be finite.
In the optimal point found by NA-IPM instead, $x_1$ and $x_2$ are infinite, which tells the user the original problem was unbounded.
It may be right to say that, in the current problem, the additional constraint introduced by Equation \eqref{eq:standard_enlargement} is equivalent to the constraint $x_1+x_2\le\alpha$ (more precisely to $x_1+x_2\le\alpha-3$), which would probably have been the first choice of anyone at the first look of Figure \ref{fig:ex_2_1}.
\begin{table}[!ht]
        \centering
        \caption{Iterations of NA-IPM solving the problem in \eqref{eq:unbounded1}}
        \begin{tabular}{|c|c|c|c|}
        \hline
        \textbf{iter} & $\bm{\mu} \in \me$ & $\bm{x} \in \me ^2$ & $\bm{f(x)} \in \me$\\
        \hline
        0 & $0.20 \alpha^{2}$ & $\begin{bmatrix}0.46 \alpha & 0.51 \alpha\end{bmatrix}$ & $0.25 \alpha^{2} - 9.67e\um1 \alpha$ \\
        \hline
        1 & $0.03 \alpha^{2}$ & $\begin{bmatrix}0.30 \alpha & 0.32 \alpha\end{bmatrix}$ & $1.55e\um3 \alpha^{2} - 6.18e\um1 \alpha$ \\
        \hline
        2 & $0.02 \alpha^{2}$ & $\begin{bmatrix}0.31 \alpha & 0.32 \alpha\end{bmatrix}$ & $2.55e\um5 \alpha^{2} - 6.35e\um1 \alpha$ \\
        \hline
        3 & $2.03e\um5 \alpha^{2}$ & $\begin{bmatrix}0.31 \alpha & 0.32 \alpha\end{bmatrix}$ & $2.55e\um7 \alpha^{2} - 6.36e\um1 \alpha$ \\
        \hline
        4 & $2.03e\um7 \alpha^{2}$ & $\begin{bmatrix}0.31 \alpha & 0.32 \alpha\end{bmatrix}$ & $-0.64 \alpha$ \\
        \hline
        5 & $7.40e\um10 \alpha^{2}$ & $\begin{bmatrix}0.31 \alpha & 0.32 \alpha\end{bmatrix}$ & $-0.64 \alpha$ \\
        \hline
        6 & $0.01 \alpha$ & $\begin{bmatrix}0.46 \alpha - 1.45  & 0.47 \alpha - 1.15 \end{bmatrix}$ & $-0.92 \alpha + 3.28 $ \\
        \hline
        7 & $2.54e\um4 \alpha $ & $\begin{bmatrix}0.49 \alpha - 1.63  & 0.50 \alpha - 1.33 \end{bmatrix}$ & $-1.00 \alpha + 3.01 $ \\
        \hline
        8 & $2.55e\um6 \alpha$ & $\begin{bmatrix}0.49 \alpha - 1.65  & 0.51 \alpha - 1.35 \end{bmatrix}$ & $-1.00 \alpha + 3.00 $ \\
        \hline
        9 & $2.55e\um8 \alpha$ & $\begin{bmatrix}0.49 \alpha - 1.65  & 0.51 \alpha - 1.35 \end{bmatrix}$ & $-1.00 \alpha + 3.00 $ \\
        \hline
        10 & $1.99e\um9$ & $\begin{bmatrix}0.49 \alpha - 1.65  & 0.51 \alpha - 1.35 \end{bmatrix}$ & $-1.00 \alpha + 3.00 $ \\
        \hline
        \end{tabular}
\label{tab:unbounded_iterations}
\end{table}

The other side of the medal is the problem described in Equation \eqref{eq:unbounded2} and drawn in Figure \ref{fig:ex_2_2}.
In this case, the primal problem is infeasible, which means that now the dual is unbounded.
Leveraging Proposition \ref{prop:weights} again, the enlarged problems becomes the one in Equation \eqref{eq:unbounded2_embed}.
Running NA-IPM in this extended problem, one appreciates that $x_5$ is equal to $\alpha$ in the optimal solution, i.e., it is nonzero and the original primal problem is infeasible.
As before, another indicator of the primal problem infeasibility is the magnitude of $\lambda_1$ and $\lambda_2$, which testifies the dual problem unboundedness.
In particular, they are equal to $\lambda_1=-0.18\alpha+0.27,\, \lambda_2=-0.09\alpha-0.36$.
The algorithm iterations are not reported for brevity.\\
\begin{minipage}{.5\linewidth}
    \vspace{3mm}
    \begin{maxi}|s|
    {x}{\quad\; x_1+x_2}{\label{eq:unbounded2}}{}
    \addConstraint{2x_1-x_2+x_3}{=-2}{}
    \addConstraint{-x_1+2x_2+x_4}{=-1}{}
    \addConstraint{x}{\ge0}{}
    \end{maxi}
\end{minipage}
\begin{minipage}{.5\linewidth}
    \vspace{3mm}
    \begin{maxi}|s|
    {x}{\quad\; x_1+x_2-\alpha x_5}{\label{eq:unbounded2_embed}}{}
    \addConstraint{2x_1-x_2+x_3-4x_5}{=-2}{}
    \addConstraint{-x_1+2x_2+x_4-3x_5}{=-1}{}
    \addConstraint{-x_3-x_4-x_6}{=-\alpha}{}
    \addConstraint{x}{\ge0}{}
    \end{maxi}
\end{minipage}

All in all, applying the embedding in \eqref{eq:standard_enlargement} the existence of optimal primal-dual solutions is guaranteed and an implementation of NA-IPM which does not check repeatedly for infeasibility can be used, helping a lot the performance.

\subsection{Experiment 3: two-objective QP}
\label{sec:ex_3}
\begin{figure*}[ht]
\begin{minipage}{.48\linewidth}
    \centering
    \includegraphics[width=\linewidth]{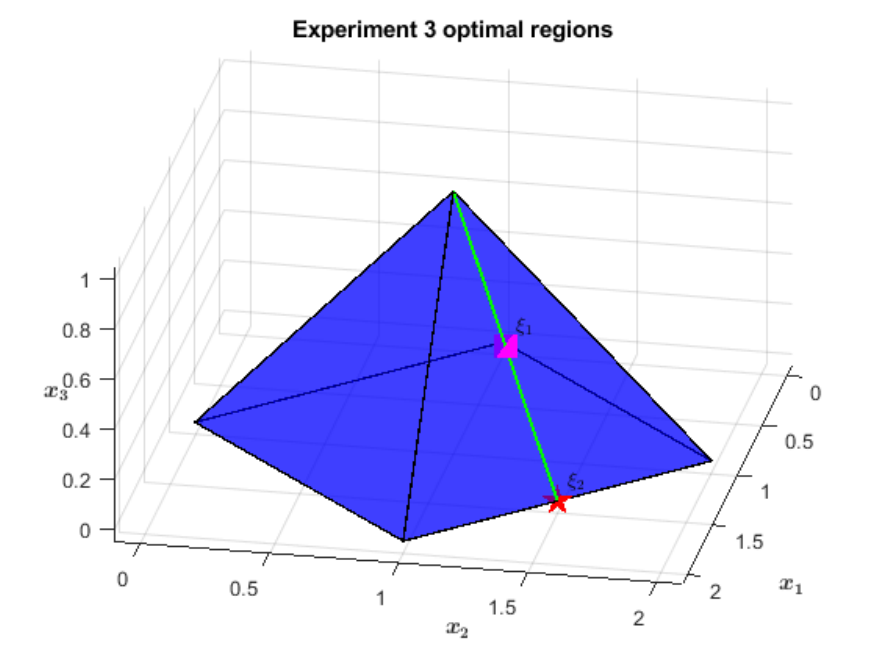}
    \caption{The segment in green is the optimal region for the first objective and $\xi_1$ is its middle point. The starred one, $\xi_2$, is the global optimum instead.}
    \label{fig:ex_3_1}
\end{minipage}
\begin{minipage}{.04\linewidth}\end{minipage}
\begin{minipage}{.48\linewidth}
    \centering
    \includegraphics[width=\linewidth]{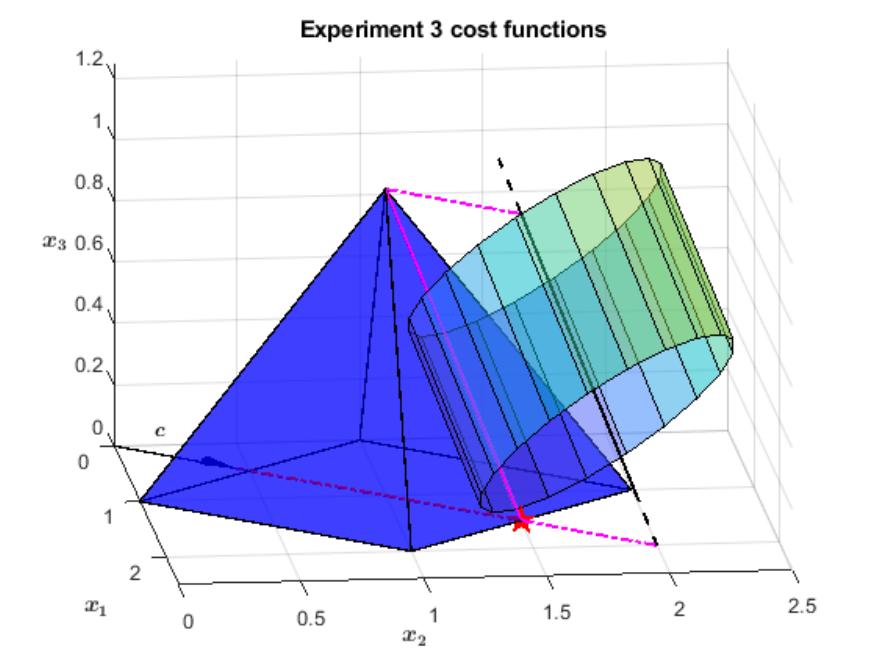}
    \caption{Cost functions: the primary is the distance from the oblique line, each infinite cylinder is a level surface; the secondary maximizes the sum of $x_1$ and $x_2$.}
    \label{fig:ex_3_2}
\end{minipage}
\end{figure*}

The problem faced in this subsection is nonlinear in the first objective and linear in the second one.
Its formal construction is reported in Appendix \ref{app:cost_ex_3}, while its standard description is in Equation \eqref{eq:2obj_QP_std} and its non-Archimedean one is in Equation \eqref{eq:2obj_QP}:\\
\vspace{2mm}
\begin{minipage}{.5\linewidth}
    \vspace{3mm}
\begin{mini}|s|
    {x}{\frac{1}{2}x^T Q x + q^T x,\; c^T x}{\label{eq:2obj_QP_std}}{\text{\hspace{-1cm}lex}}
    \addConstraint{-x_1+x_2+x_3}{\le 1}{}
    \addConstraint{-x_1-x_2+x_3}{\le 1}{}
    \addConstraint{x_1-x_2+x_3}{\le 1}{}
    \addConstraint{x_1+x_2+x_3}{\le 3}{}
    \addConstraint{x_3}{\ge 0}{}
\end{mini}
\end{minipage}
\begin{minipage}{.5\linewidth}
    \vspace{3mm}
    \begin{mini}|s|
    {x}{\frac{1}{2}x^T Q\eta x + (c+q\eta)^Tx}{\label{eq:2obj_QP}}{}
    \addConstraint{-x_1+x_2+x_3}{\le 1}{}
    \addConstraint{-x_1-x_2+x_3}{\le 1}{}
    \addConstraint{x_1-x_2+x_3}{\le 1}{}
    \addConstraint{x_1+x_2+x_3}{\le 3}{}
    \addConstraint{x_3}{\ge 0}{}
    \end{mini}
\end{minipage}
\begin{equation*}
    Q = 
    \begin{bmatrix}
    \;\,10  & -2       & 4 \\
    -2      & \;\,10   & 4 \\
    \;\;\,4 &  \;\;\,4 & 4 
    \end{bmatrix},
    \quad\quad
    q = 
    \begin{bmatrix}
    -16\\-16\\-16
    \end{bmatrix},
    \quad \quad
    c =
    \begin{bmatrix}
    -1\\-1\\\;\;\,0
    \end{bmatrix}.
\end{equation*}
The primal feasible region is a right square pyramid of height one, see Figure \ref{fig:ex_3_1}.
Its basis lies on the $x_1x_2$ plane, has center in $(1,\,1,\,0)$ and vertices in $(1,\,0,\,0)$, $(0,\,1,\,0)$, $(1,\,2,\,0)$, $(2,\,1,\,0)$.
Moving to the cost function, the first objective penalizes points far from the line $\mathcal{L}=(t,\,t,\,-2t+4)$, $t\in\mr$, i.e., the axis of the cylinder in Figure \ref{fig:ex_3_2}.
The optimal region, highlighted in green
in both the figures, is the height of the furthest-from-origin face of the pyramid, identified by the intersection of the latter with the cylinder of radius $\tfrac{1}{\sqrt{3}}$ and $\mathcal{L}$ as axis.
Finally, the second objective selects from the whole apothem the point as close as possible to the $x_1x_2$ plane, namely $\xi_2=(\tfrac{1}{2},\,\tfrac{1}{2},\,0)$.

Before reaching the global optimum, NA-IPM passes close to the midpoint of the first objective optimal region.
In the current problem it is $\xi_1=(\tfrac{5}{4},\,\tfrac{5}{4},\,\tfrac{1}{2})$, which is approached since iteration 3, as can be seen in Table \ref{tab:2obj_QP_iterations}.
Then, starting from iteration 6, it goes towards $\xi_2$ since the algorithm perceives the first objective as fully optimized, as one can deduce from the change in magnitude of $\mu$.
\begin{table}[!ht]
        \centering
        \caption{Iterations of NA-IPM solving the problem in \eqref{eq:2obj_QP}}
        \begin{tabular}{|c|c|c|c|}
        \hline
        \textbf{iter} & $\bm{\mu} \in \mr$ & $\bm{x} \in \mr^3$ & $\bm{f(x)} \in \me$\\
        \hline
        0 & $3.37 $ & $\begin{bmatrix}1.46 & 1.46 & 1.46 \end{bmatrix}$ & $-31.72  - 2.92 \eta$ \\
        \hline
        1 & $1.23 $ & $\begin{bmatrix}1.26 & 1.26 & 0.66 \end{bmatrix}$ & $-30.64  - 2.52 \eta$ \\
        \hline
        2 & $0.11 $ & $\begin{bmatrix}1.27 & 1.27 & 0.39 \end{bmatrix}$ & $-29.71  - 2.54 \eta$ \\
        \hline
        3 & $2.12e\um3 $ & $\begin{bmatrix}1.30 & 1.30 & 0.39 \end{bmatrix}$ & $-29.99  - 2.60 \eta$ \\
        \hline
        4 & $2.12e\um5 $ & $\begin{bmatrix}1.30 & 1.30 & 0.39 \end{bmatrix}$ & $-30.00  - 2.61 \eta$ \\
        \hline
        \textbf{5} & $\mathbf{2.12e\um7} $ & $\begin{bmatrix}\bm{1.30} & \bm{1.30} & \bm{0.39} \end{bmatrix}$ & $\bm{-30.00  - 2.61 \eta}$ \\
        \hline
        \textbf{6} & $\bm{0.15 \eta}$ & $\begin{bmatrix}\mathbf{1.38} & \mathbf{1.38} & \mathbf{0.24} \end{bmatrix}$ & $\bm{-30.00  - 2.76 \eta}$ \\
        \hline
        7 & $0.06 \eta$ & $\begin{bmatrix}1.50 & 1.50 & 0.00 \end{bmatrix}$ & $-30.00  - 3.00 \eta$ \\
        \hline
        8 & $1.14e\um3 \eta$ & $\begin{bmatrix}1.50 & 1.50 & 0.00 \end{bmatrix}$ & $-30.00  - 3.00 \eta$ \\
        \hline
        9 & $1.14e\um5 \eta$ & $\begin{bmatrix}1.50 & 1.50 & 0.00 \end{bmatrix}$ & $-30.00  - 3.00 \eta$ \\
        \hline
        10 & $1.14e\um7 \eta$ & $\begin{bmatrix}1.50 & 1.50 & 0.00 \end{bmatrix}$ & $-30.00  - 3.00 \eta$ \\
        \hline
    \end{tabular}
\label{tab:2obj_QP_iterations}
\end{table}

\subsection{Experiment 4: three-objective QP}
\label{sec:ex_4}
\begin{figure*}[ht]
\begin{minipage}{.48\linewidth}
    \centering
    \includegraphics[width=\linewidth]{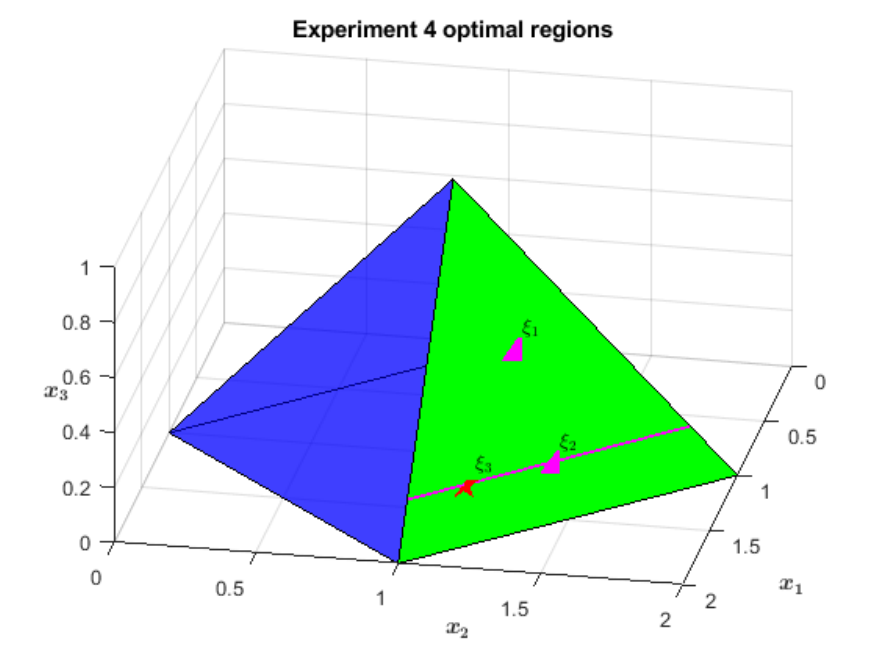}
    \caption{The surface in green is the optimal region of the first objective ($\xi_1$ is its midpoint); the segment in magenta is the optimal region also for the second objective ($\xi_2$ is its midpoint); $\xi_3$ is the global optimum.}
    \label{fig:ex_4_1}
\end{minipage}
\begin{minipage}{.04\linewidth}
\end{minipage}
\begin{minipage}{.48\linewidth}
    \centering
    \includegraphics[width=\linewidth]{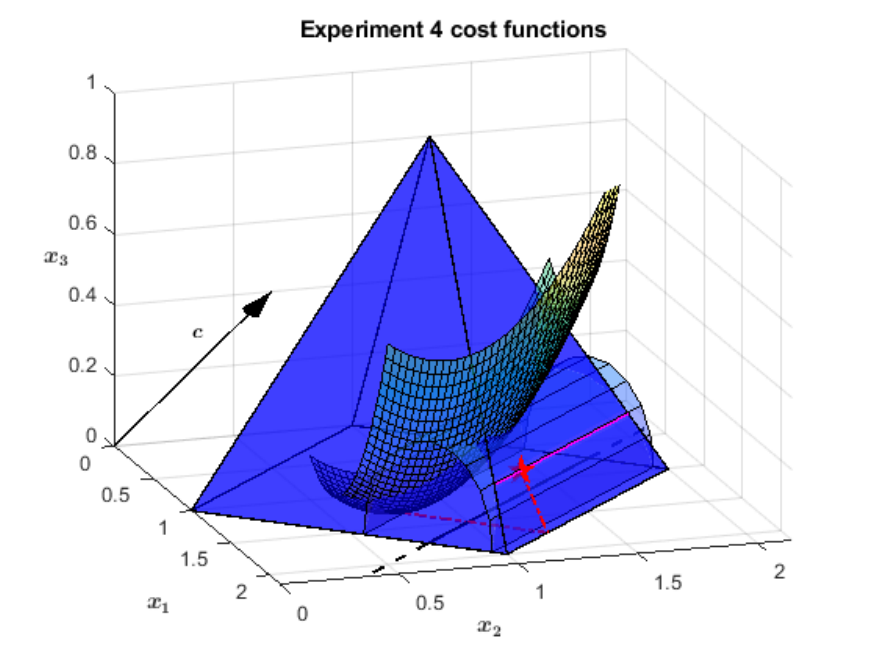}
    \caption{Cost functions: the primary objective is the vector $c$; the secondary one is the distance form the black line on the $x_1x_2$ axis (the cylinder is a level surface); the tertiary one is the infinite paraboloid centered in $(\tfrac{5}{4},\,\tfrac{3}{4},\,0)$.}
    \label{fig:ex_4_2}
\end{minipage}
\end{figure*}

The last test problem involves three distinct objective functions.
Again, the step by step construction is reported in appendix (Appendix \ref{app:cost_ex_4} to be precise), while its standard and non-Archimedean descriptions are the ones in Equation \eqref{eq:3obj_QP} and \eqref{eq:NA_3obj_QP}, respectively:\\
\vspace{2mm}
\begin{minipage}{.5\linewidth}
    \vspace{3mm}
    \begin{mini}|s|
{x}{c^T x,\; \frac{1}{2}x^T Q x + q^T x,\; \frac{1}{2} x^T P x + p^T x}{\label{eq:3obj_QP}}{\text{\hspace{-1cm}lex}}
    \addConstraint{-x_1+x_2+x_3}{\le 1}{}
    \addConstraint{-x_1-x_2+x_3}{\le 1}{}
    \addConstraint{x_1-x_2+x_3}{\le 1}{}
    \addConstraint{x_1+x_2+x_3}{\le 3}{}
    \addConstraint{x_3}{\ge 0}{}
\end{mini}
\end{minipage}
\begin{minipage}{.5\linewidth}
    \vspace{3mm}
    \begin{mini}|s|
{x}{\frac{1}{2}x^T (Q\eta+P\eta^2) x + (c+q\eta+ p\eta^2)^T x}{\label{eq:NA_3obj_QP}}{}
    \addConstraint{-x_1+x_2+x_3}{\le 1}{}
    \addConstraint{-x_1-x_2+x_3}{\le 1}{}
    \addConstraint{x_1-x_2+x_3}{\le 1}{}
    \addConstraint{x_1+x_2+x_3}{\le 3}{}
    \addConstraint{x_3}{\ge 0}{}
\end{mini}
\end{minipage}
where
\begin{equation}
    c = 
    \begin{bmatrix}
    -1\\-1\\-1
    \end{bmatrix},
    \quad
    Q = 
    \begin{bmatrix}
    2 & 2 & 0\\
    2 & 2 & 0\\
    0 & 0 & 4
    \end{bmatrix},
    \quad
    q =
    \begin{bmatrix}
    -5 \\ -5 \\ \;\;\,0
    \end{bmatrix},
    \quad 
    P = 
    \begin{bmatrix}
    4 & 0 & 0\\
    0 & 4 & 0\\
    0 & 0 & 0
    \end{bmatrix},
    \quad 
    p =
    \begin{bmatrix}
    -5 \\ -3 \\ \;\;\,2
    \end{bmatrix}.
\label{eq:3obj_QP_values}
\end{equation}
The feasible region is the same pyramid as in the previous problem, Section \ref{sec:ex_3}.
This time however, there are three objectives and only the first one is linear (see Figure \ref{fig:ex_4_2}).
The latter promotes as optimal the whole farthest-from-origin face of the pyramid, the green one in Figure \ref{fig:ex_4_1}.
The second objective is again a penalization for the distance from a line, which in this case is $\mathcal{L}=(t,\,-t+\tfrac{5}{2},\,0)$, $t\in\mr$.
Of the original triangle, the only surface optimal also for the second objective is the segment linking $(\tfrac{23}{12},\,\tfrac{11}{12},\,\tfrac{1}{6})$ and $(\tfrac{11}{12},\,\tfrac{23}{12},\,\tfrac{1}{6})$, highlighted in magenta in both the figures.
The third objective selects just one point of that segment: $\xi_3=(\tfrac{5}{3},\,\tfrac{7}{6},\,\tfrac{1}{6})$.
It consists in a potential function with the shape of an infinite paraboloid, centered in $(\tfrac{5}{4},\,\tfrac{3}{4},\,0)$ and growing linearly along $\vv{x_3}$.

As one can immagine, NA-IPM passes close to two more points before reaching $\xi_3$.
The first of those is the center of the green triangle, i.e., $\xi_1=(\tfrac{5}{4},\,\tfrac{5}{4},\,\tfrac{1}{2})$.
Then, it goes towards $\xi_2=(\tfrac{8}{3},\,\tfrac{8}{3},\,\tfrac{1}{6})$, which would have been one of the global optima (but the unique one with maximum centrality) if the problem had had just two objectives rather than three.
Indeed, it is the midpoint of the second objective optimal segment along the green triangle.
Table \ref{tab:3obj_QP_iterations} identifies this switch in the optimizing direction at iterations 5-6, when $\mu$ becomes infinitesimal of the first order.
Eventually, the third objective takes the stage and brings the optimization away from $\xi_2$, pointing in the direction of $\xi_3$.
This happens at iterations 10-11, when the magnitude of $\mu$ decreases again, making it infinitesimal of the second order.
Such a value of $\mu$ testifies the completion of both the first and second objective optimization, leaving space to the third.
\begin{table}[!ht]
        \centering
        \caption{Iterations of NA-IPM solving the problem in \eqref{eq:3obj_QP}}
        \begin{tabular}{|c|c|c|c|}
        \hline
        \textbf{iter} & $\bm{\mu} \in \mr$ & $\bm{x} \in \mr^3$ & $\bm{f(x)} \in \me$\\
        \hline
        0 & $0.53 $ & $\begin{bmatrix}1.46 & 1.46 & 1.46 \end{bmatrix}$ & $-4.38  - 3.64 \eta - 6.08 \eta^2$ \\
        \hline
        1 & $0.21 $ & $\begin{bmatrix}1.32 & 1.32 & 0.74 \end{bmatrix}$ & $-3.37  - 10.27 \eta - 5.08 \eta^2$ \\
        \hline
        2 & $0.02 $ & $\begin{bmatrix}1.30 & 1.30 & 0.40 \end{bmatrix}$ & $-3.01  - 11.83 \eta - 4.44 \eta^2$ \\
        \hline
        3 & $1.60e\um4 $ & $\begin{bmatrix}1.30 & 1.30 & 0.40 \end{bmatrix}$ & $-3.00  - 11.84 \eta - 4.44 \eta^2$ \\
        \hline
        4 & $1.60e\um6 $ & $\begin{bmatrix}1.30 & 1.30 & 0.40 \end{bmatrix}$ & $-3.00  - 11.84 \eta - 4.44 \eta^2$ \\
        \hline
        \textbf{5} & $\mathbf{1.61e\um8} $ & $\begin{bmatrix}\bm{1.30} & \bm{1.30} & \bm{0.40} \end{bmatrix}$ & $\bm{-3.00  - 11.84 \eta - 4.44 \eta^2}$ \\
        \hline
        \textbf{6} & $\bm{0.06 \eta}$ & $\begin{bmatrix}\bm{1.38} & \bm{1.38} & \bm{0.25} \end{bmatrix}$ & $\bm{-3.00  - 12.13 \eta - 3.92 \eta^2}$ \\
        \hline
        7 & $2.21e\um3 \eta$ & $\begin{bmatrix}1.41 & 1.41 & 0.17 \end{bmatrix}$ & $-3.00  - 12.17 \eta - 3.66 \eta^2$ \\
        \hline
        8 & $2.46e\um5 \eta $ & $\begin{bmatrix}1.42 & 1.42 & 0.17 \end{bmatrix}$ & $-3.00  - 12.17 \eta - 3.64 \eta^2$ \\
        \hline
        9 & $2.48e\um7 \eta $ & $\begin{bmatrix}1.42 & 1.42 & 0.17 \end{bmatrix}$ & $-3.00  - 12.17 \eta - 3.64 \eta^2$ \\
        \hline
        \textbf{10} & $\bm{1.62e}\mathbf{\um}\bm{9 \eta}$ & $\begin{bmatrix}\bm{1.42} & \bm{1.42} & \bm{0.17} \end{bmatrix}$ & $\bm{-3.00  - 12.17 \eta - 3.64 \eta^2}$ \\
        \hline
        \textbf{11} & $\bm{0.14 \eta^2}$ & $\begin{bmatrix}\bm{1.54} & \bm{1.29} & \bm{0.17} \end{bmatrix}$ & $\bm{-3.00  - 12.17 \eta - 3.82 \eta^2}$ \\
        \hline
        12 & $0.01 \eta^2$ & $\begin{bmatrix}1.65 & 1.19 & 0.17 \end{bmatrix}$ & $-3.00  - 12.17 \eta - 3.89 \eta^2$ \\
        \hline
        13 & $1.63e\um4 \eta^2$ & $\begin{bmatrix}1.67 & 1.17 & 0.17 \end{bmatrix}$ & $-3.00  - 12.17 \eta - 3.89 \eta^2$ \\
        \hline
        14 & $1.78e\um6 \eta^2$ & $\begin{bmatrix}1.67 & 1.17 & 0.17 \end{bmatrix}$ & $-3.00  - 12.17 \eta - 3.89 \eta^2$ \\
        \hline
        15 & $1.59e\um8 \eta^2$ & $\begin{bmatrix}1.67 & 1.17 & 0.17 \end{bmatrix}$ & $-3.00  - 12.17 \eta - 3.89 \eta^2$ \\
        \hline
    \end{tabular}
\label{tab:3obj_QP_iterations}
\end{table}

\subsection{An informal comparison with the scalarization approach}
\hl{A typical standard approach to tackle (deterministic) lexicographic optimization problems is the scalarization 
method} \cite{Natural_computing}\hl{.
In terms of effectiveness, it is not able to guarantee equivalence between the original problem and the scalarized one, even if it converges to a solution significantly faster.
More important, the choice of the scalarization parameters always requires a trial-and-error process to achieve a good tuning: even in case of a weighted combination of the objectives, the motto ``the smaller the better'' can be applied.
Table} \ref{tab:comparison_scalar} \hl{testifies this fact quantitatively, using the benchmark in Section} \ref{sec:ex_4} \hl{as reference problem.}
\begin{table}[ht]
    \centering
    \caption{\hl{The approximation error of the optimal lexicographic point $\xi_3$ does not monotonically decrease when reducing the scalarization weight $w$ on problem} \eqref{eq:3obj_QP}.}
    \begin{tabular}{c|c|c|c|c|c}
        $\mathbf{w}$    & 1e-1 & 1e-2 & 1e-3 & 1e-4 & 1e-5 \\
                        \hline
        \textbf{Error}  & 6.59e-3 & 6.78e-4 & 8.17e-5 & 1.46e-1 & 3.29e-1\\

    \end{tabular}
    \label{tab:comparison_scalar}
\end{table}
\hl{The study replaces the non-Archimedean weight $\eta$ in} \eqref{eq:NA_3obj_QP} \hl{with a standard one ($w$), and solves the problem using the Matlab \texttt{quadrprog} routine.
Evidences highlight how too large or too small values of $w$ introduce a bias in the optimization which negatively affect the computations.
Actually, the error made by the solving algorithm , i.e., the norm of the discrepancy between the theoretically optimal point $\xi_3$ and the solution found is not acceptable even when $w$ decreases too much. 
Furthermore, the scalarization approach is always problem-dependent, which means it cannot constitute a general-purpose solving paradigm.}


\section{Conclusions}
\label{sec:conclusions}
This paper builds upon a theoretical and practical achievement happened the last years concerning the solution of lexicographic multi-objective linear programming problem.
Indeed, in \cite{CococcioniEtAlAMC2018}, it has been demonstrated that a non-Archimedean version of the simplex algorithm is able to solve that problem, in an elegant and powerful way.
The aim of this work was to verify whether or not a non-Archimedean version of an interior point algorithm (called by the authors NA-IPM) is able to deliver a solution to the same class of problems and go beyond it.
This paper answers affirmatively, and this is interesting for situations where an interior point method is known to perform better than the simplex, like in high dimensional problems.
With this NA-IPM at hand, it was natural to verify if it was able to solve a convex lexicographic multi-objective quadratic programming problem, and again the answer is positive.

The proposed implementation of NA-IPM shows a polynomial complexity, and its implementation guarantees to converge in finite time.
The new algorithm also enjoys a very light embedding which gets rid of the issues related to infeasibility and unboundedness.
As a practical application, this paper considered and discussed lexicographic optimization problems as well as infeasible/unbounded ones, testing NA-IPM on four linear and quadratic programming examples, achieving the expected results.
The successful solution of such problems
paves the way to more difficult ones,
such as lexicographic multi-objective semi-definite programming problems, which is left as a future study. 
It should also be noted that the NA-IPM can be used to solve other family of problems involving infinitesimal/infinite numbers,
as testified by Section \ref{sec:ex_2}.
The study of its performances on harder examples
is left for a future study as well.

\appendix

\setcounter{section}{0}
\setcounter{equation}{0}
\section[\appendixname~\thesection]{Cost function construction in experiment 3}
\label{app:cost_ex_3}
This appendix shows the construction of the cost function in problem \eqref{eq:2obj_QP}.
The second objective discussion is omitted since it is linear and therefore trivial.

The idea is to make the apothem of the farthest-from-origin face of the pyramid optimal for the first objective.
This segment belongs to the line $\mathcal{G}=(t,\,t,\,-2y+3)$, $t\in\mr$.
A possible cost function is the one which penalizes the distance from a line $\mathcal{L}$ parallel to $\mathcal{G}$.
In addition, $\mathcal{L}$ must satisfy the property that no point in the primal feasible region (the pyramid) is closer to it than $\mathcal{G}$.
A feasible choice is the translation of $\mathcal{G}$ along $\vv{x_3}$, say $\mathcal{L} = \mathcal{G}+[0,\,0,\,1]= (t,\,t,\,-2t+4)$, $t\in\mr$.
In this case, the line versor is $v=(-\tfrac{1}{\sqrt{6}},\, -\tfrac{1}{\sqrt{6}},\, \tfrac{\sqrt{2}}{\sqrt{3}})$.

The distance $d$ of any point $x\in\mr^3$ from $\mathcal{L}$ comes from the following equation
\begin{equation}
    d = x^1 - (v^Tx^1)v,
\label{eq:distance}    
\end{equation}
where $x^1=x-x^0$ and $x^0$ is any point on $\mathcal{L}$, i.e., $x^0\in\mathcal{L}$.
Few lines to explain the equation follow.
The translation $x-x^0$ brings the system origin to $x^0$, which allows one to consistently execute the projection of the point of interest (identified by $x^1$ in the new reference system) along $\mathcal{L}$.
To do this, first one computes the norm of the projection by the inner product $v^T x^1$, then one constructs the projection multiplying such norm by the line reference direction, that is its versor $v$.
Since the projection of $x^1$ on $\mathcal{L}$ is parallel to it by definition, the difference between $x^1$ and such a projection must be orthogonal to it, i.e., it is the distance of $x$ from $\mathcal{L}$ in the coordinate system with origin in $x^0$.

In problem \eqref{eq:2obj_QP}, $x^0$ is arbitrarily chosen as the vector $[2,\,2,\,0]^T$; therefore, the distance is:
\begin{equation*}
    d =
    \begin{bmatrix}
    x_1-2\\x_2-2\\x_3
    \end{bmatrix} -
    \left(
    -\frac{1}{\sqrt{6}}(x_1+x_2-4)+\frac{\sqrt{3}}{\sqrt{3}}x_3
    \right)
    \begin{bmatrix}
    -\frac{1}{\sqrt{6}}\\-\frac{1}{\sqrt{6}}\\\frac{\sqrt{2}}{\sqrt{3}}
    \end{bmatrix} = 
    \frac{1}{6}
    \begin{bmatrix}
    5x_1+x_2+2x_3-8\\
    -x_1+5x_2+2x_3-8\\
    2x_1+2x_2+2x_3-8
    \end{bmatrix}.
\end{equation*}
To model the problem as a QP task, one can use the squared norm of the distance $d$ as cost function.
After having done the calculations, its expression is (up to a multiplicative factor)
\begin{equation}
    \|d\|^2 = 30x_1^2 + 30x_2^2 + 12x_3^2 - 12x_1x_2 + 24x_1x_3 + 24x_2x_3 - 96x_1 - 96x_2 - 96x_3 + k,
\label{eq:squared_distance_ex3}
\end{equation}
where $k$ is a constant with no impact on the optimization and therefore it shall not be considered furthermore.
Dividing \eqref{eq:squared_distance_ex3} by 6 one gets exactly the matrix $Q$ and the vector $q$ in \eqref{eq:2obj_QP} (remind that $Q$ is divided by 2 in the cost function and therefore the diagonal entries must be doubled).

The optimal solution of the whole problem $\overline{x}$ must lie on the pyramid face apothem due to the first objective function, i.e., $\overline{x}=(\overline{t},\,\overline{t},\,-2\overline{t}+3)$ for some $\overline{t}\in[1,\,\tfrac{3}{2}]$.
The second objective, namely $c$ in \eqref{eq:2obj_QP}, forces $\overline{x}$ (read $\overline{t}$) to posses another property: it must also solve the following scalar optimization problem t:
\begin{maxi*}|s|
    {t}{2t}{}{}
    \addConstraint{t\in\left[1,\,\frac{3}{2}\right]}{}{}
\end{maxi*}
whose optimal solution is $\overline{t}=\tfrac{3}{2}$.
Therefore, the optimal solution of the whole problem is $\overline{x}=[\tfrac{3}{2},\,\tfrac{3}{2},\,0]^T$, which coincides with $\xi_2$ in Section \ref{sec:ex_3}.

\section[\appendixname~\thesection]{Cost function construction in experiment 4}
\label{app:cost_ex_4}

This appendix shows the construction of the cost function in problem \eqref{eq:3obj_QP}.
It consists of three different objectives, the first is linear while the other two are quadratic.
Let one disclose them in order.

The first objective indicates as optimal surface the whole farthest-from-origin face of the pyramidal feasible region.
To select it, the cost vector $c$ must be orthogonal to the plane $\pi$ the triangular face belongs to.
Actually, that plane is identified by the equation $x_1+x_2+x_3=3$, as testified by the fourth constraint in \eqref{eq:3obj_QP}.
Basic notions in linear algebra say that the normal direction to a plane is exactly the span of the vector filled with its coefficients, i.e., $[1,\,1,\,1]^T$.
Since the problem is a minimization one and the constraint is a lesser than or equal to one, then the normal vector must be reversed, that is $c=[-1,\,-1,\,-1]^T$, exactly as in \eqref{eq:3obj_QP_values}.

Of such an optimal face, the second objective selects just a segment parallel to the basis.
To do it, the distance from a straight line is used, similarly to what done in Appendix \ref{app:cost_ex_3}.
Actually, the reference line must be parallel to both the basis and the face; $\mathcal{L}=(t,\,-t+\tfrac{5}{2},\,0)$, $t\in\mr$ is a feasible choice.
Indeed, it is parallel to the $x_1x_2$ plane (read the pyramid basis) since it is constant on $x_3$, while it is parallel to the face because it is parallel to its basis $\mathcal{B}=(t,\,-t+3,\,0)$, $t\in\mr$.
Using \eqref{eq:distance}, the distance equation in this case becomes:
\begin{equation*}
    d = 
    \begin{bmatrix}
    x_1-\frac{7}{4}\\
    x_2-\frac{3}{4}\\
    x_3
    \end{bmatrix} - \frac{1}{\sqrt{2}}
    \left(1-x_1+x_2\right)
    \begin{bmatrix}
    - \frac{1}{\sqrt{2}}\\  \frac{1}{\sqrt{2}}\\ 0
    \end{bmatrix} = \frac{1}{2}
    \begin{bmatrix}
    x_1+x_2-\frac{5}{2}\\ x_1+x_2-\frac{5}{2}\\ 2x_3
    \end{bmatrix},
\end{equation*}
where $x^0=[\frac{7}{4},\,\frac{3}{4},\,0]^T$ and $v=[-\tfrac{1}{\sqrt{2}},\, \tfrac{1}{\sqrt{2}},\,0]^T$.
The squared norm of the distance is now (up to a multiplicative factor)
\begin{equation}
    \|d\|^2 = x_1^2+x_2^2+2x_3^2+2x_1x_2-5x_1 - 5x_2 + k
    \label{eq:squared_distance_ex4}
\end{equation}
which originates $Q$ and $q$ of \eqref{eq:3obj_QP_values}.

To analytically identify the optimal region for the second objective, let one project the problem on the $x_1=x_2$ plane and then retrieve the whole surface leveraging the fact that such surface is a segment parallel to $\mathcal{L}$.
Applying the condition $x_1=x_2=t$, the pyramid face collapses to its apothem, whose equation is $\mathcal{G}=(t,\,t,\,-2t+3)$, $t\in[1,\,\tfrac{3}{2}]$ (see Appendix \ref{app:cost_ex_3}).
Since the optimal region must belong to $\pi$, to retrieve it one has to substitute the equation of $\mathcal{G}$ into \eqref{eq:squared_distance_ex4}, to differentiate with respect to $t$ and to impose the optimality condition:
\begin{equation*}
\begin{gathered}
    x\in\mathcal{G} \Longrightarrow \|d\|^2 = 12t^2-34t+k,\;t\in\left[1,\,\frac{3}{2}\right], \\
    \frac{\text{d}\|d\|^2}{\text{d}t} = 24t-34 = 0 \Longrightarrow t = \frac{17}{12}\in\left[1,\,\frac{3}{2}\right].
\end{gathered}
\end{equation*}
Therefore, the optimal surface for the second objective is a line parallel to $\mathcal{L}$ and passing through $\xi=[\tfrac{17}{12},\,\tfrac{17}{12},\,\tfrac{1}{6}]^T$ (by construction it already belongs to $\pi$).
Basic linear algebra considerations say that this request is equivalent to constrain the sought region to the line $\mathcal{S}=(t,\,-t+\tfrac{17}{6},\,\tfrac{1}{6})$.
To identify the values of $t$ for which $\mathcal{S}$ belongs to the triangular face, one needs to intersect $\mathcal{S}$ with the pyramid edges $\mathcal{E}_0=(t,\,1,\,-t+2)$ and $\mathcal{E}_1=(1,\,t,\,-t+2)$.
The intersecting points locate at $t=\tfrac{11}{6}$ and $t=1$, respectively.
Therefore the optimal region for the second objective is $\mathcal{S}=(t,\,-t+\tfrac{17}{6},\,\tfrac{1}{6})$, $t\in[1,\,\tfrac{11}{6}]$, henceforth indicated by $\overline{\mathcal{S}}$.

The third objective is a potential function.
By arbitrary choice, it is an infinite paraboloid centered in $[\tfrac{5}{4},\,\tfrac{3}{4},\,0]^T$ and growing linearly along $\vv{x_3}$, i.e., its equation is (scaled up for practical reasons by a factor 2)
\begin{equation}
    h(x) = 2\left(\left(x_1-\frac{5}{4}\right)^2 + \left(x_2-\frac{3}{4}\right)^2 + x_3\right) = 
    2x_1^2+2x_2^2-5x_1-3x_2+2x_3 + k
    \label{eq:potential}
\end{equation}
which coincides with $P$ and $p$ of \eqref{eq:3obj_QP_values}.

The optimal point $\overline{x}$ of the whole problem is the point of $\overline{\mathcal{S}}$ (optimal surface identified by the second objective) which minimizes \eqref{eq:potential}.
Since $\overline{x}\in\overline{\mathcal{S}}$, it has the form $\overline{x}=(\overline{t},\, -\overline{t}+\tfrac{17}{6},\, \tfrac{1}{6})$, $\overline{t}\in[1,\,\tfrac{11}{6}]$.
Substituting this parametric description into \eqref{eq:potential} and applying the first order optimality condition, one gets
\begin{equation*}
    x\in\overline{\mathcal{S}} \Longrightarrow h(t) = 4t^2-\frac{40}{3}t+k,\;t\in[1,\,\tfrac{11}{6}] \quad \quad \quad \frac{\text{d}h}{\text{d}t} = 8t-\frac{40}{3} = 0 \Longrightarrow \overline{t} = \frac{5}{3}\in[1,\,\tfrac{11}{6}].
\end{equation*}
Therefore, $\overline{x}=[\tfrac{5}{3},\,\tfrac{7}{6},\,\tfrac{1}{6}]^T$.

\section*{Acknowledgements}
The authors wish to express our warmest appreciation to prof. Vieri Benci for his generous help in teaching us Robinson's non-Standard Analysis and his Alpha-Theory.

This work has been funded by the Italian Ministry of Education and Research (MIUR) in the framework of the CrossLab project (Departments of Excellence).

The authors are very thankful to the three anonymous reviewers for their accurate comments and insightful remarks.

\bibliographystyle{plain}
\bibliography{OperationsResearch.bib,Grossone.bib,GameTheory.bib,NonArchimedeanAnalysis.bib}

\end{document}